\documentclass{article}
\usepackage{maa-monthly2}
\usepackage{mathtools,xspace}
\usepackage{tikz}
\usetikzlibrary{arrows,positioning,automata}
\usepackage[outline]{contour}
\contourlength{1.2pt}
\usepackage{framed}
\usepackage{colonequals}
\usepackage{calc}
\usepackage[clockwise]{rotating}
\usepackage{fancyhdr}
\usepackage{multirow}
\usepackage[export]{adjustbox}
\usepackage{floatpag}
\usepackage{wrapfig}
\usepackage{floatrow}
\usepackage{microtype}
\usepackage[pdftitle={The looping rate and sandpile density of planar graphs},
 pdfauthor={Adrien Kassel and David B. Wilson},
 pdfpagescrop={37 38 540.2 756},
 bookmarks=true,bookmarksopen=true,bookmarksopenlevel=2]{hyperref}
\usepackage{bookmark}

\newcommand{\Z}{\mathbb{Z}}
\newcommand{\Q}{\mathbb{Q}}
\newcommand{\E}{\mathbb{E}}

\newcommand{\G}{\mathcal{G}}
\newcommand{\arcsec}{\operatorname{arcsec}}
\newcommand{\level}{\operatorname{level}}

\begin{document}

\title{The looping rate and sandpile density\\ of planar graphs}
\markright{\quad\quad\quad\quad\quad\quad\quad Looping rate and sandpile density of planar graphs}
\author{Adrien Kassel and David B.~\!Wilson}
\date{}

\maketitle
\begin{abstract}
  We give a simple formula for the looping rate of
  loop-erased random walk on a finite planar graph.  The looping rate is closely related
  to the expected amount of sand in a recurrent sandpile on the graph.
  The looping rate formula is well-suited to taking limits
  where the graph tends to an infinite lattice, and we use it to give an
  elementary derivation of the (previously computed) looping rate and
  sandpile densities of the square, triangular, and honeycomb lattices,
  and compute (for the first time) the looping rate and
  sandpile densities of many other lattices, such as the kagom\'e lattice,
  the dice lattice, and the truncated hexagonal lattice
  (for which the values are all rational), and the square-octagon lattice
  (for which it is transcendental).
\end{abstract}

\section{Spanning trees and sandpiles}\label{intro}

\begin{wrapfigure}{r}{0.32\textwidth}
  \begin{center}
    \vspace*{-2.5\baselineskip}
    \includegraphics[width=\textwidth]{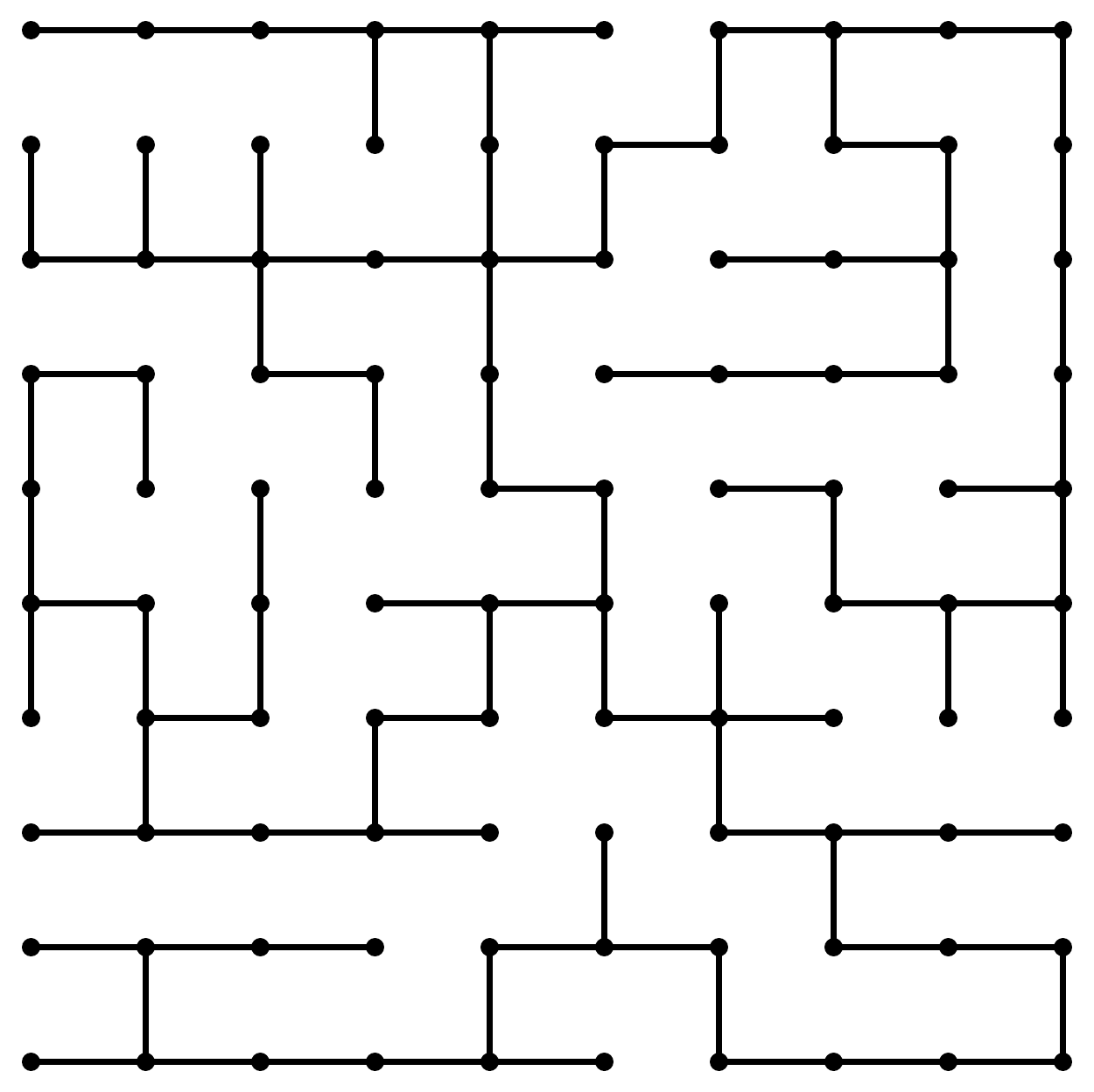}
  \end{center}
\small{{\sc Figure~1.} Uniformly random spanning tree (UST) of a $10\times 10$ grid.}
\label{fig:UST}
\end{wrapfigure}
\addtocounter{figure}{1}
Spanning trees on graphs have a long history which goes back to
Kirchhoff, who used them to compute effective resistances in electric
networks \cite{Kirchhoff}.  Formally, a spanning tree of a finite connected graph
is a collection of edges through which any two vertices may be
connected, and which contains no cycles.  As we shall explain
later, certain electrical quantities in a resistive electric network
correspond to the probabilities of certain events in a uniformly
random spanning tree (UST).

\begin{figure}[b!]
\begin{center}
\includegraphics[width=\textwidth]{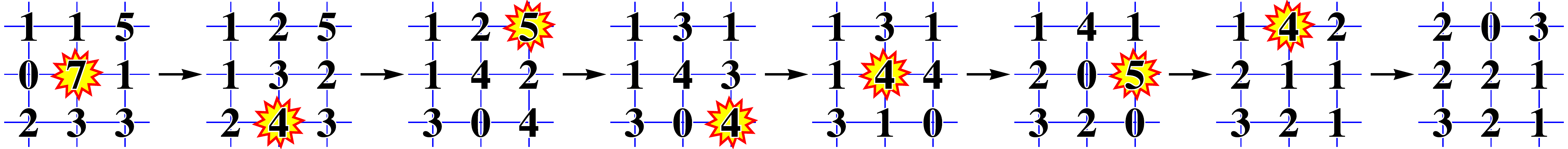}
\end{center}
\caption{Stabilization of a sandpile configuration on a $3\times 3$ grid,
where the sink (not shown) is the region outside the grid.}
\end{figure}
Uniformly random spanning trees are also closely related to the abelian
sandpile model of self-organized criticality \cite{BTK}, as was shown
by Dhar and Majumdar \cite{MR1044086,Majumdar-Dhar}.  In the abelian
sandpile model on a finite graph, each vertex has a non-negative
integer number of grains of sand.  If the vertex contains at least as
many grains of sand as it has neighbors, then the vertex is
\textit{unstable}, and may \textit{topple}, sending one grain of sand to
each neighbor.  Usually there is a designated sink vertex which never
topples.  Assuming every vertex is connected to the sink, then we may
repeatedly topple unstable vertices
until every vertex is stable.  The
resulting sandpile is called the stabilization of the original
sandpile, and is independent of the order in which vertices are
toppled (which is the abelian property).
Some sandpile configurations
are \textit{recurrent}, meaning that from any sandpile configuration, it
is possible to add some amount of sand to the vertices and stabilize
to obtain the given configuration.  These sandpile configurations are
the recurrent states of the Markov chain which at each step adds a grain
of sand to a random vertex and then stabilizes the configuration.
Majumdar and Dhar gave a bijection between the recurrent sandpile
configurations of a finite graph $\G$ with given sink~$s$ and the
spanning trees of~$\G$ \cite{majumdar-dhar:height}, which we will discuss further in
Section~\ref{sec:sandpile}.

Pemantle~\cite{Pemantle} initiated the study of uniformly random spanning
trees on the infinite lattice $\Z^d$.  Of course there are infinitely
many such spanning trees, so some care is needed to make sense of
this.  Pemantle considered a sequence of finite graphs $(\G_n)_{n\ge 1}$ which
converges to $\Z^d$, and argued that the distribution of uniform
spanning trees on~$\G_n$ converges in a suitable sense, and took the
limit to be the definition of the uniform spanning tree on $\Z^d$.  We
say that the sequence~$(\G_n)_{n\ge 1}$ converges to $\Z^d$ if for every finite induced
subgraph $H$ of $\Z^d$, for sufficiently large~$n$ we have that $H$ is
contained in~$\G_n$ as an induced subgraph.  For any finite box
$B_m=\{-m,\mbox{$-m+1$},\dots,\mbox{$m-1$},m\}^d$ centered at the origin,
and for those $n$'s
that are sufficiently large for $B_m\subset\G_n$, we can consider a
uniformly random spanning tree $T_n$ on $\G_n$ restricted to the box~$B_m$.
The restriction $T_n|_{B_m}$ naturally contains no cycles, but need
not be connected.  Pemantle showed that the distribution of the set of
edges in the restricted tree $T_n|_{B_m}$ converges as $n\to\infty$, and
that this limiting distribution is independent of the choice of
sequence $(\G_n)_{n\ge 1}$ converging to $\Z^d$.  Since there is a canonical
limiting distribution on acyclic sets of edges within each box
centered at the origin, taken together they define a random forest on
$\Z^d$, which is called the uniform spanning forest USF$(\Z^d)$.
Pemantle showed that for $d\leq 4$, with probability~$1$ the USF
contains just a single tree, in which case it is called the uniform
spanning tree UST$(\Z^d)$, but that for $d\geq 5$, with probability~$1$
the USF contains infinitely many trees.  Each such tree has
a path leading to infinity, and one point of view is that the trees
are connected through infinity.
(See \cite{BLPS} and \cite{LP:book} for further developments and
streamlined proofs.)

\begin{figure}[b]
\floatbox[{\capbeside\thisfloatsetup{capbesideposition={right,bottom},capbesidewidth=4.5cm}}]{figure}[\FBwidth]
{\caption{A portion of a uniformly random spanning tree on $\Z^2$.
  The restriction of the spanning tree to the box is a forest,
  with each connected component reaching the boundary.
(This configuration was produced using an algorithm which computes probabilities of local spanning tree events on the infinite lattice.)
}}
{\hspace{0pt}\includegraphics[width=0.5\textwidth]{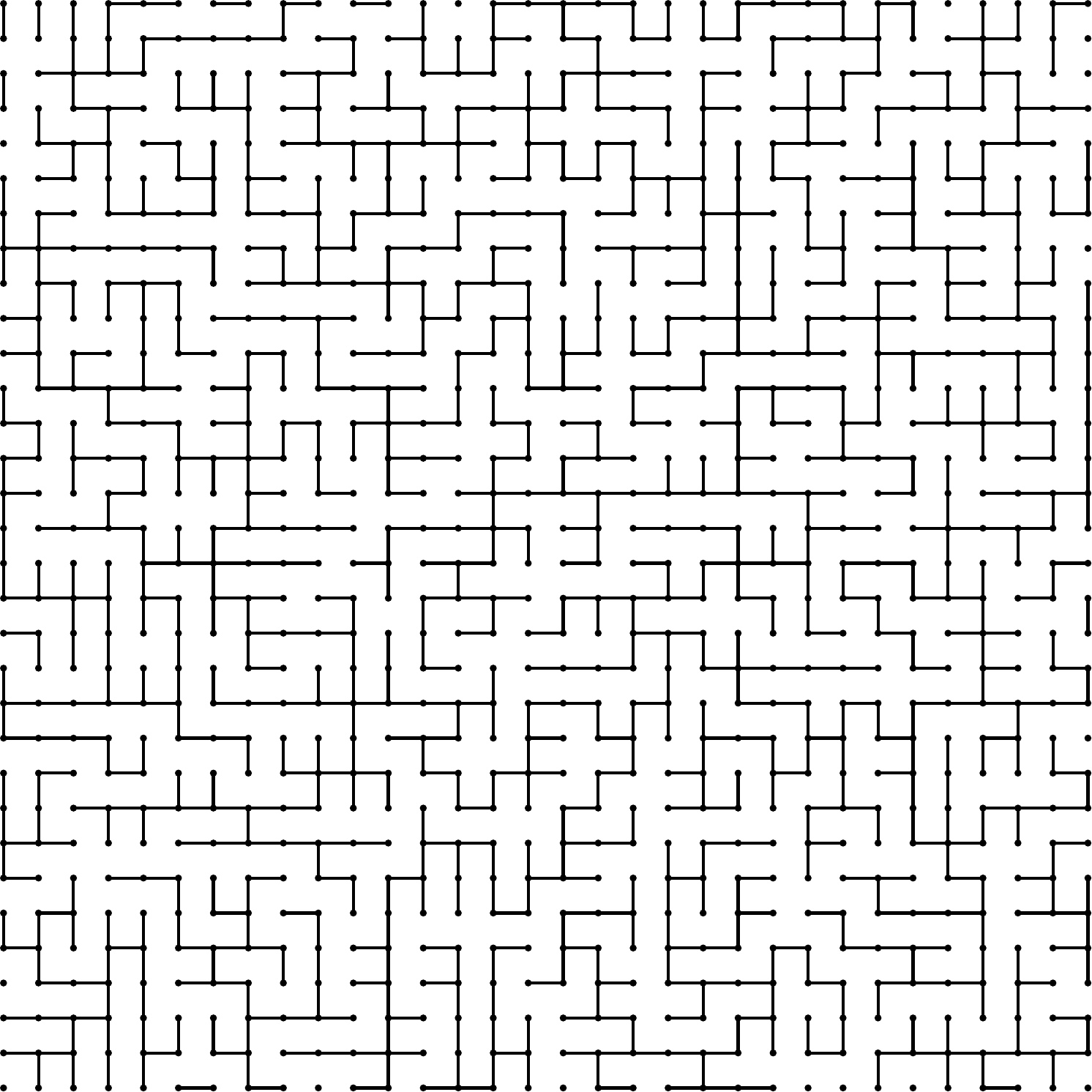}\hspace{20pt}}
\end{figure}

Burton and Pemantle \cite{BP} showed how to compute, for any finite
collection of edges, the probability that USF$(\Z^d)$ contains
those edges.  These probabilities can be expressed in terms of the
discrete Green's function, and for UST$(\Z^2)$, they are all rational
polynomials in $1/\pi$.

Using their bijection between spanning trees and sandpiles, and the
ability to compute local probabilities for spanning trees, Majumdar
and Dhar \cite{majumdar-dhar:height} showed that the probability that a
vertex of $\Z^2$ has zero grains of sand is $2/\pi^2-4/\pi^3$.
Computing the other sandpile height probabilities turned out to be
much harder.  The reason is that the maps between spanning trees and
sandpiles are not local in nature, so that local events for spanning
trees do not correspond to local events for sandpiles, except in some
special cases.  Nonetheless, the sandpile height probabilities were
computed by Priezzhev \cite{priezzhev:heights,Priezzhev}, although the
expressions he gave contained a singular integral involving
trigonometric functions.  Grassberger (in unpublished work) evaluated
these integrals numerically, and made the surprising observation that
for $\Z^2$, the average amount of sand per vertex, called the
\textit{sandpile density}, was numerically indistinguishable from $17/8$.
Despite much effort, it took eighteen years for this $17/8$
conjecture to be proved \cite{JPR,PPR:54,KW4} (see also
\cite{caracciolo-sportiello:heights}).

While the $17/8$ conjecture was fully proved, none of the three
proofs really explained why the answer was rational, since they
went through
calculations with intermediate values involving $1/\pi$ or
integrals of trigonometric functions, and when these intermediate
values were combined to give the final answer, the transcendental
parts mysteriously cancelled.

\begin{sidewaystable}[htbp]
\newcommand{\depends}{\multirow{2}{*}{\raisebox{8pt}{\small \parbox{\widthof{start vertex}}{\centering depends on \\[-2pt] start vertex}}}}
\newcommand{\lattice}[1]{\multirow{2}{*}{\includegraphics[height=33pt,margin=0pt -3.5pt]{#1-small}}}
\newcommand{\rowgap}{-8pt}
\newcommand{\ngap}{6pt}
\begin{center}
\resizebox{\textheight}{!}{%
\begin{tabular}{|c@{\,}c@{\,}|c|c|c|c|c|c|}
\hline
\multicolumn{2}{|c|}{}&&&&&&\\[-11pt]
  \multicolumn{2}{|c|}{\multirow{2}{*}{lattice}} &
  \parbox{\widthof{\raisebox{0pt}{$\displaystyle\frac{\text{unicycle}}{\text{tree}\times\text{edge}}$} ratio}}{\raisebox{0pt}{$\displaystyle\frac{\text{unicycle}}{\text{tree}\times\text{edge}}$} ratio} &
  \parbox{\widthof{\!mean unicycle\!}}{\centering \makebox[0pt][c]{mean unicycle}\\[-2pt] loop length} &
  \parbox{\widthof{mean LERW}}{\centering mean LERW\\[-2pt] loop length} &
  \parbox{\widthof{LERW looping rate}}{\centering discrete-time\\LERW looping rate} &
  \parbox{\widthof{UST path to $\infty$}}{\centering mean number\\[-2pt]of neighbors on\\[-2pt] UST path to $\infty$} &
  \parbox{\widthof{sandpile density}}{\centering sandpile density\\[-2pt]} \\[12pt]
  \multicolumn{2}{|c|}{} &
  $\tau$ & $\lambda=\Pr[e\notin T]/\tau$ & $1/\rho$ & $\rho=\tau+\frac12\Pr[e\in T]$ & $\delta\rho$ & $\bar\sigma=(\delta\rho+\delta-1)/2$ \\[2pt]
\hline
\hline
&&&&&&&\\[\rowgap]
  \multirow{2}{*}{square} &
  \lattice{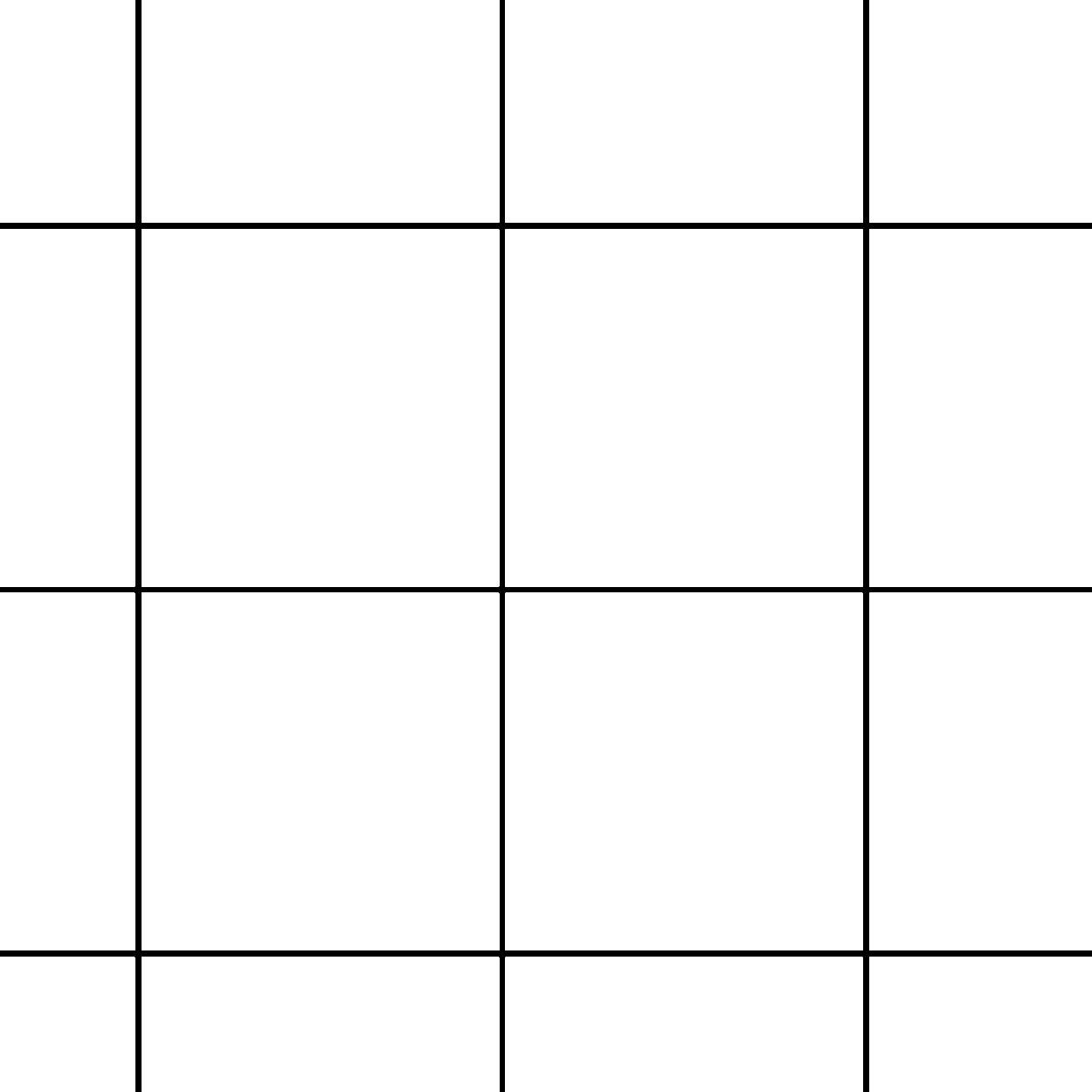} &
  $1/16$ &
  \multirow{2}{*}{$8$} &
  $16/5$ &
  $5/16$ &
  $5/4$ &
  $17/8$ \\[\ngap]
  & &
  \small $0.0625\phantom{00\dots}$ & &
  \small $3.2\phantom{00000\dots}$ &
  \small $0.3125\phantom{00\dots}$ &
  \small $1.25\phantom{0000\dots}$ &
  \small $2.125\phantom{000\dots}$ \\
\hline
&&&&&&&\\[\rowgap]
  \multirow{2}{*}{triangular} &
  \lattice{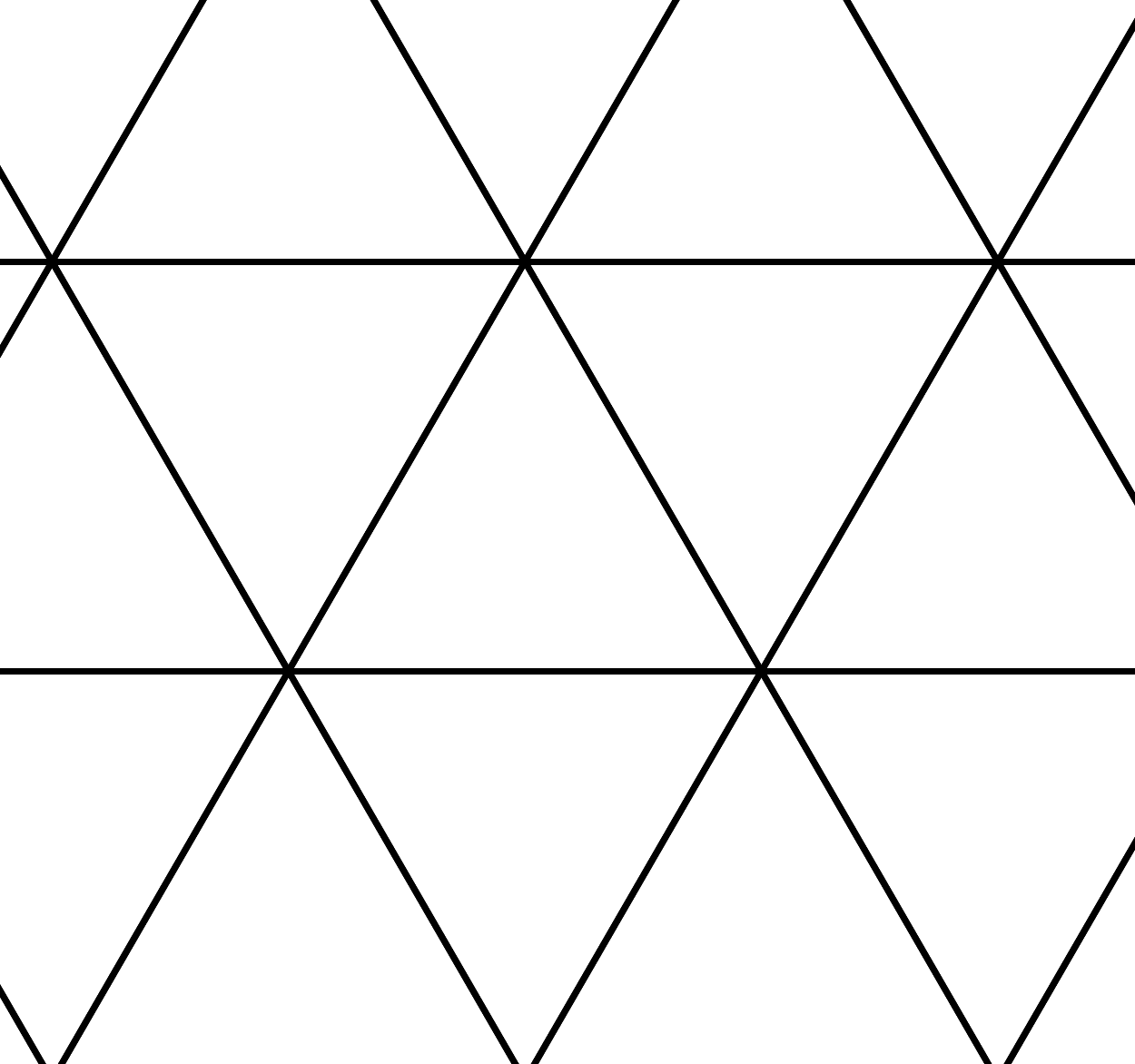} &
  $1/9$ &
  \multirow{2}{*}{$6$} &
  $18/5$ &
  $5/18$ &
  $5/3$ &
  $10/3$ \\[\ngap]
 & &
  \small $0.111111\dots$ & &
  \small $3.6\phantom{00000\dots}$ &
  \small $0.277778\dots$ &
  \small $1.666667\dots$ &
  \small $3.333333\dots$ \\
\hline
&&&&&&&\\[\rowgap]
  \multirow{2}{*}{honeycomb} &
  \lattice{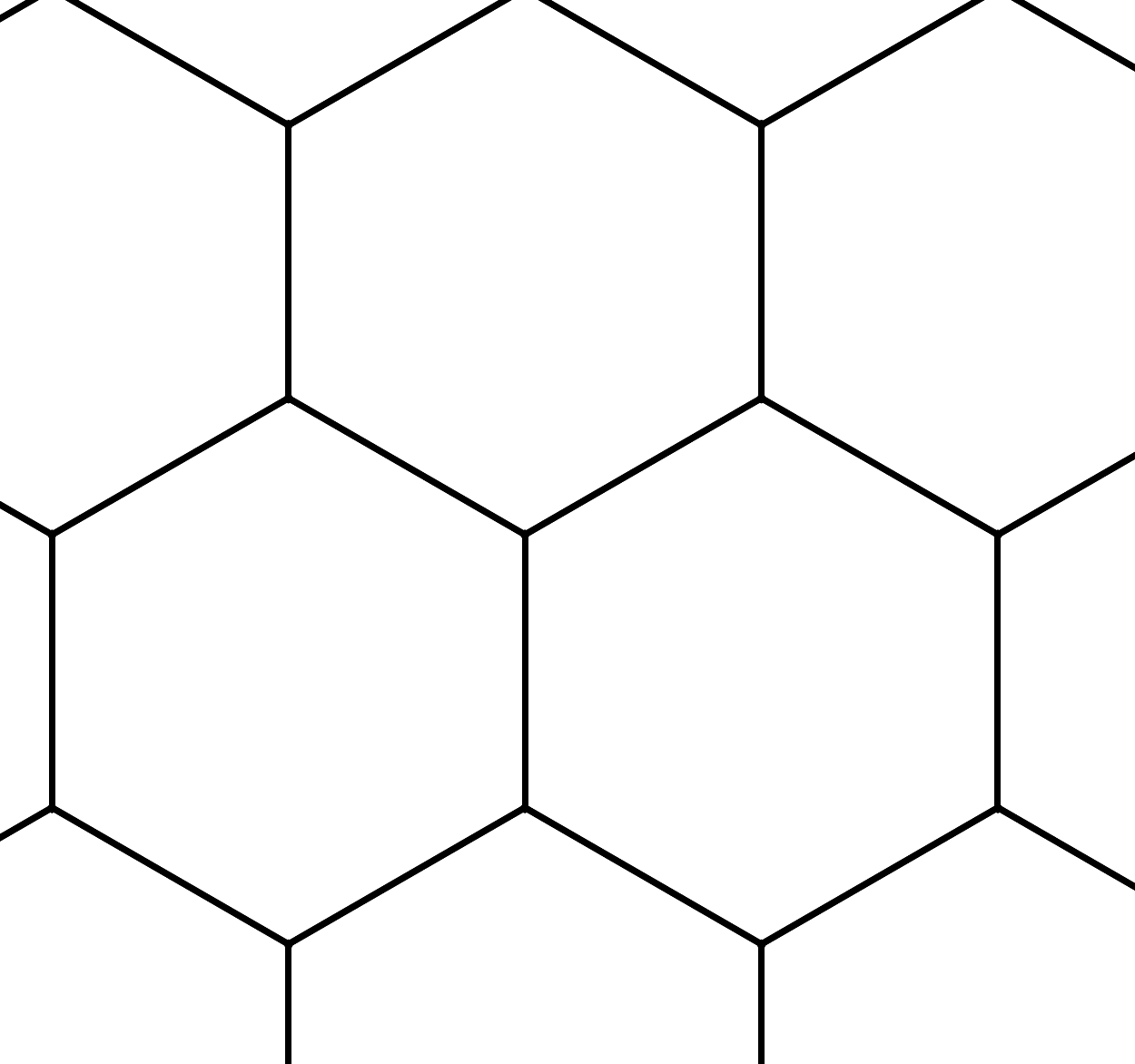} &
  $1/36$ &
  \multirow{2}{*}{$12$} &
  $36/13$ &
  $13/36$ &
  $13/12$ &
  $37/24$ \\[\ngap]
 & &
  \small $0.027778\dots$ & &
  \small $2.769231\dots$ &
  \small $0.361111\dots$ &
  \small $1.083333\dots$ &
  \small $1.541667\dots$ \\
\hline
&&&&&&&\\[\rowgap]
  \multirow{2}{*}{\parbox{\widthof{trihexagonal}}{\centering kagom\'e / \\trihexagonal}} &
  \lattice{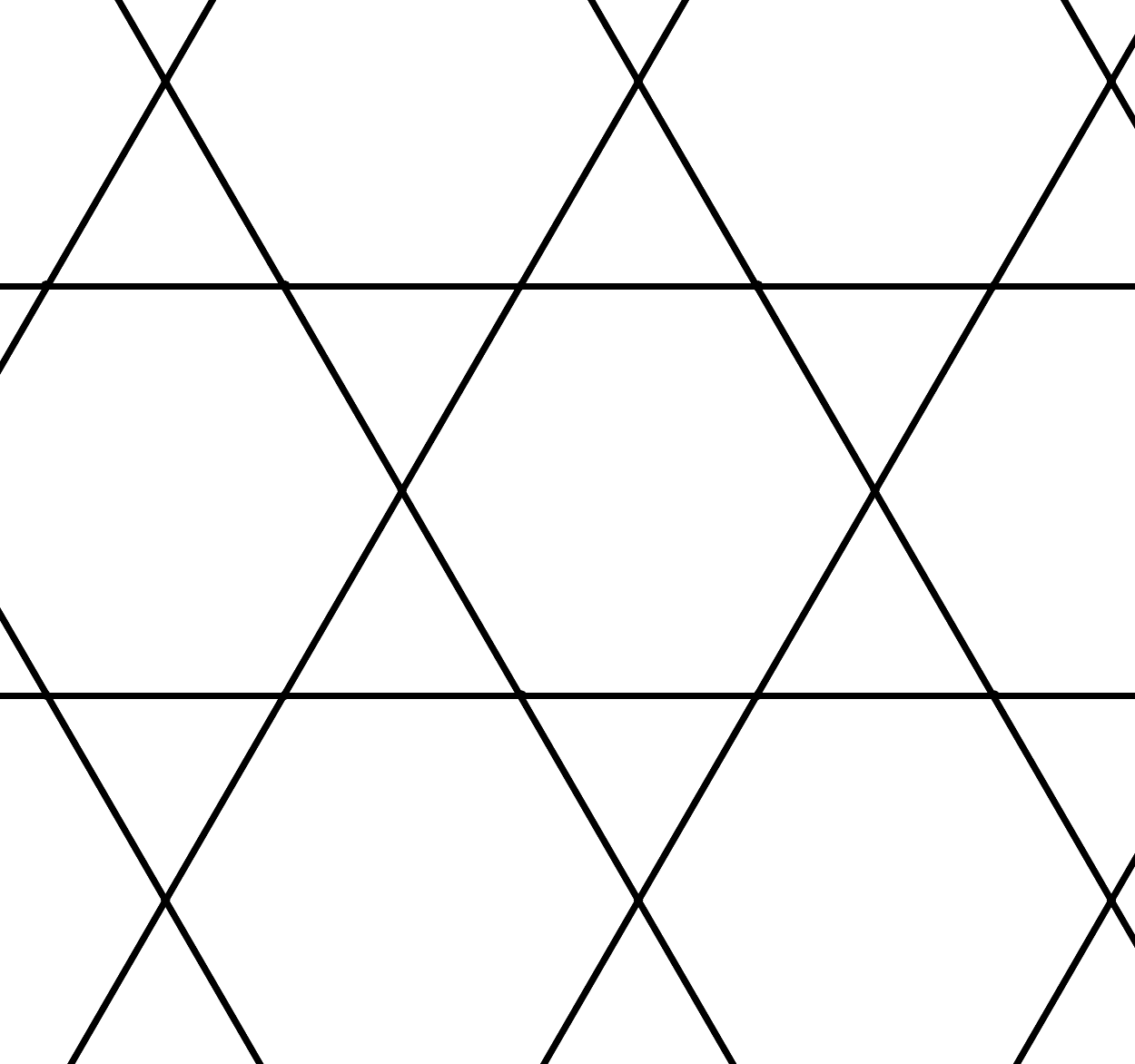} &
  $1/12$ &
  \multirow{2}{*}{$6$} &
  \multirow{2}{*}{$3$} &
  $1/3$ &
  $4/3$ &
  $13/6$ \\[\ngap]
 & &
  \small $0.083333\dots$  & & &
  \small $0.333333\dots$ &
  \small $1.333333\dots$ &
  \small $2.166667\dots$ \\
\hline
&&&&&&&\\[\rowgap]
  \multirow{2}{*}{\parbox{\widthof{rhombille}}{\centering dice / \\rhombille}} &
  \lattice{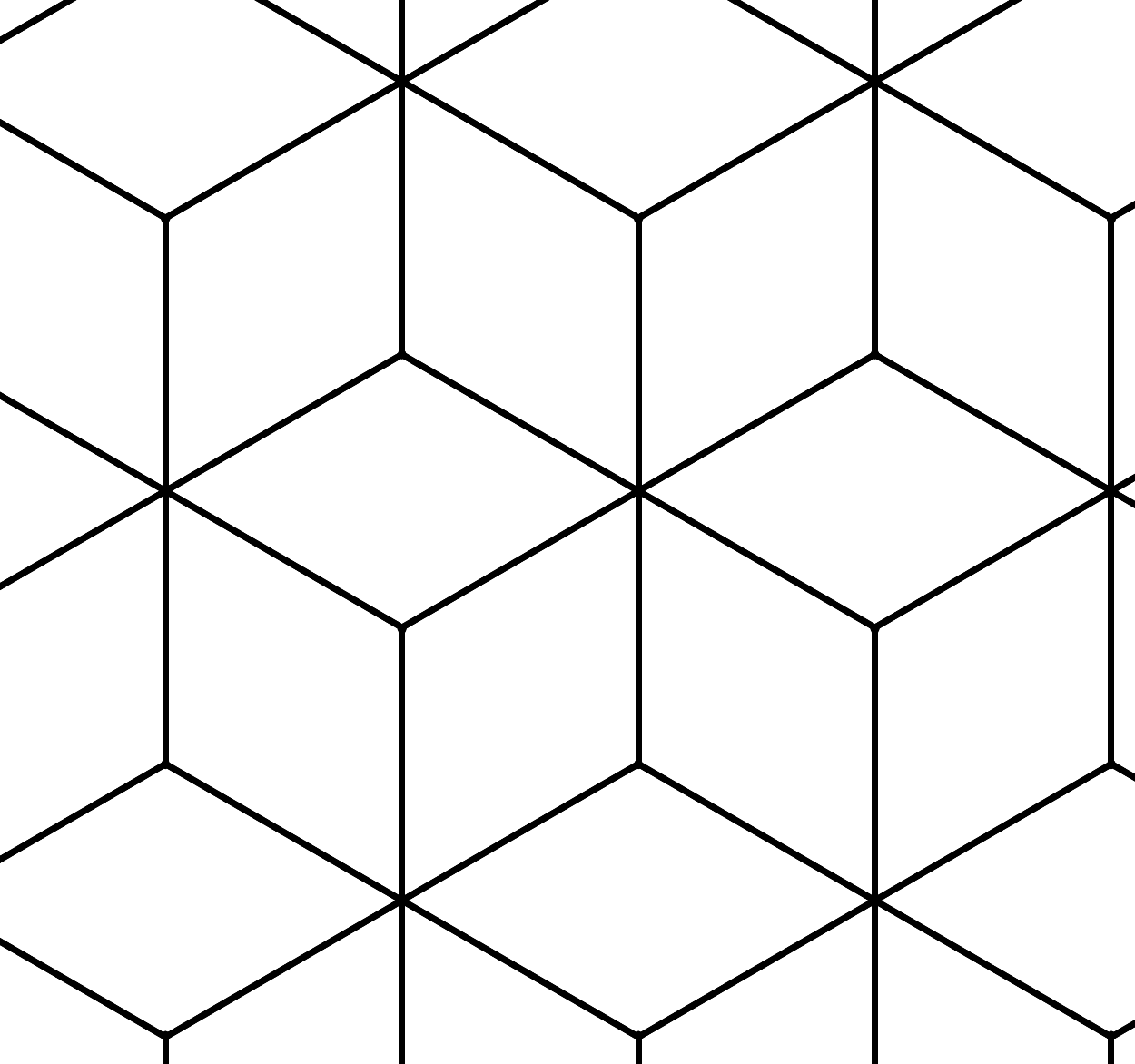} &
  $1/16$ &
  \multirow{2}{*}{$8$} &
  $16/5$ &
  $5/16$ &
  $5/4$ &
  $17/8$ \\[\ngap]
  & &
  \small $0.0625\phantom{00\dots}$ & &
  \small $3.2\phantom{00000\dots}$ &
  \small $0.3125\phantom{00\dots}$ &
  \small $1.25\phantom{0000\dots}$ &
  \small $2.125\phantom{000\dots}$ \\
\hline
&&&&&&&\\[\rowgap]
  \multirow{2}{*}{\parbox{\widthof{\small truncated hexagonal}}{\centering Fisher / \\ \small truncated hexagonal}} &
  \lattice{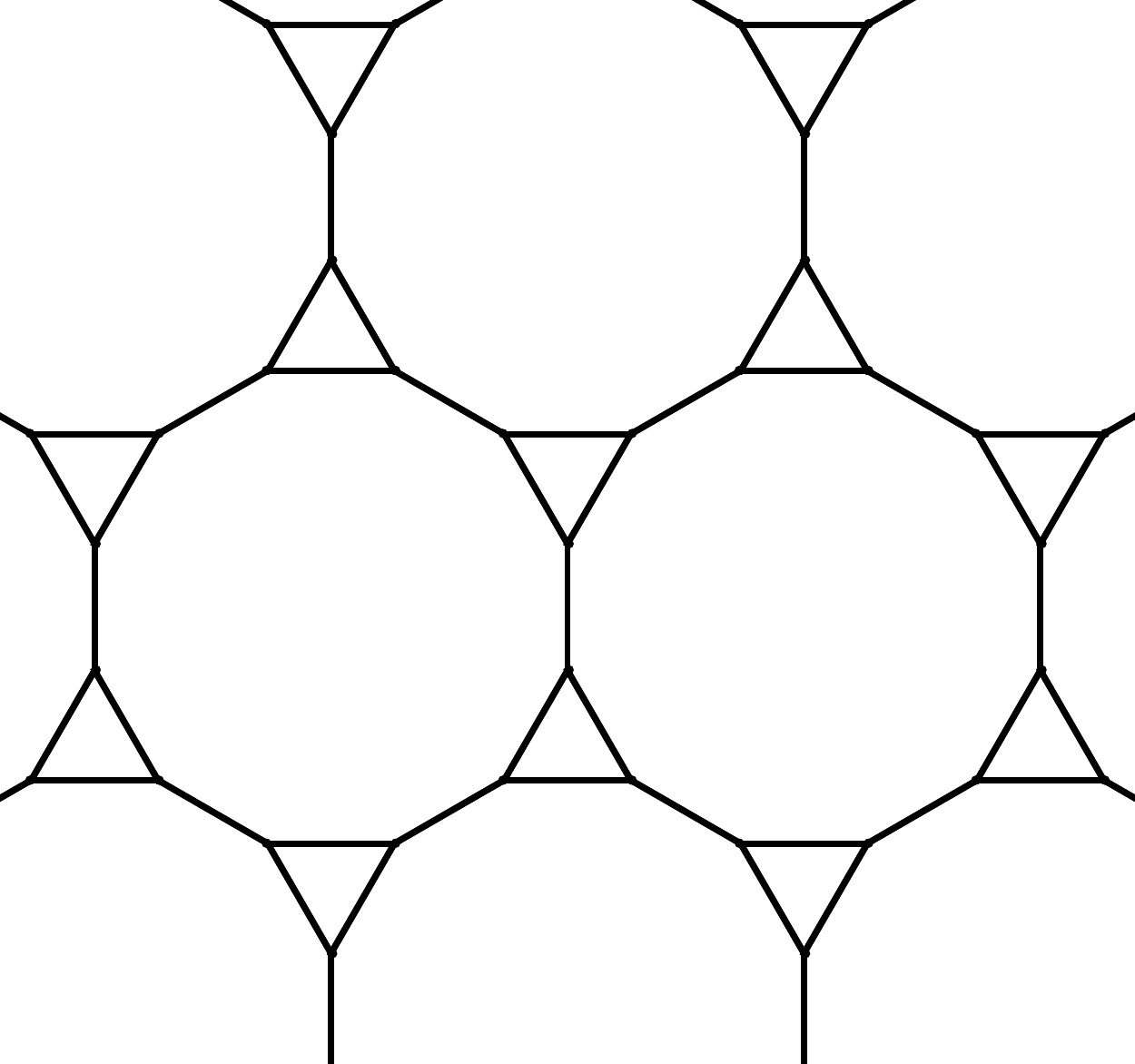} &
  $59/900$ &
  $300/59$ &
  $900/359$ &
  $359/900$ &
  $359/300$ &
  $959/600$ \\[\ngap]
 & &
  \small $0.065556\dots$ &
  \small $5.084746\dots$ &
  \small $2.506964\dots$ &
  \small $0.398889\dots$ &
  \small $1.196667\dots$ &
  \small $1.598333\dots$ \\
\hline
&&&&&&&\\[\rowgap]
\multirow{2}{*}{triakis triangular} &
  \lattice{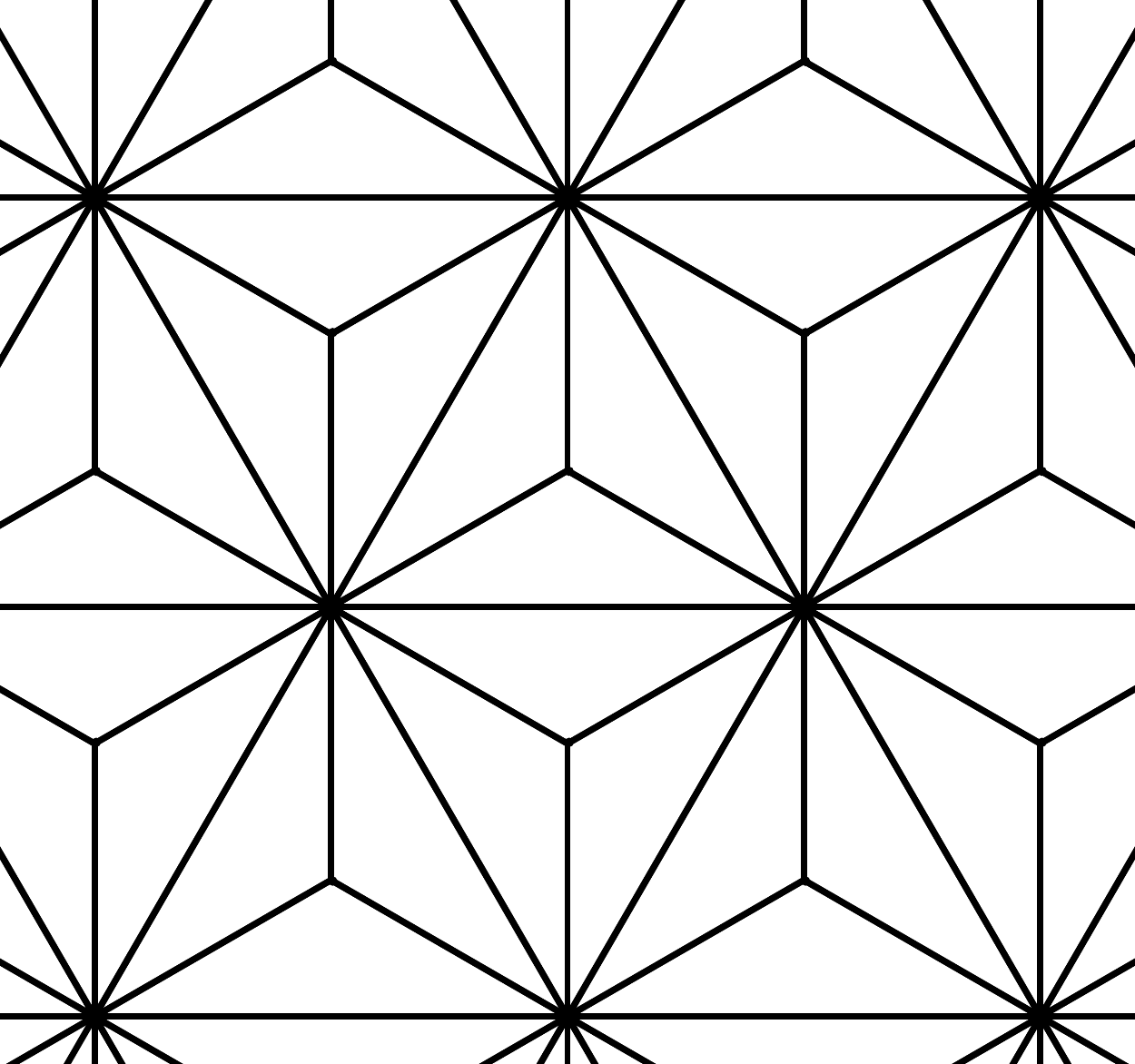} &
  $17/150$ &
  $100/17$ &
  $25/7$ &
  $7/25$ &
  $42/25$ &
  $167/50$ \\[\ngap]
 & &
  \small $0.113333\dots$ &
  \small $5.882353\dots$ &
  \small $3.571429\dots$ &
  \small $0.28\phantom{0000\dots}$ &
  \small $1.68\phantom{0000\dots}$ &
  \small $3.34\phantom{0000\dots}$ \\
\hline
&&&&&&&\\[\rowgap]
  \multirow{2}{*}{\parbox{\widthof{square-octagon /}}{\centering square-octagon / \\truncated square}} &
  \lattice{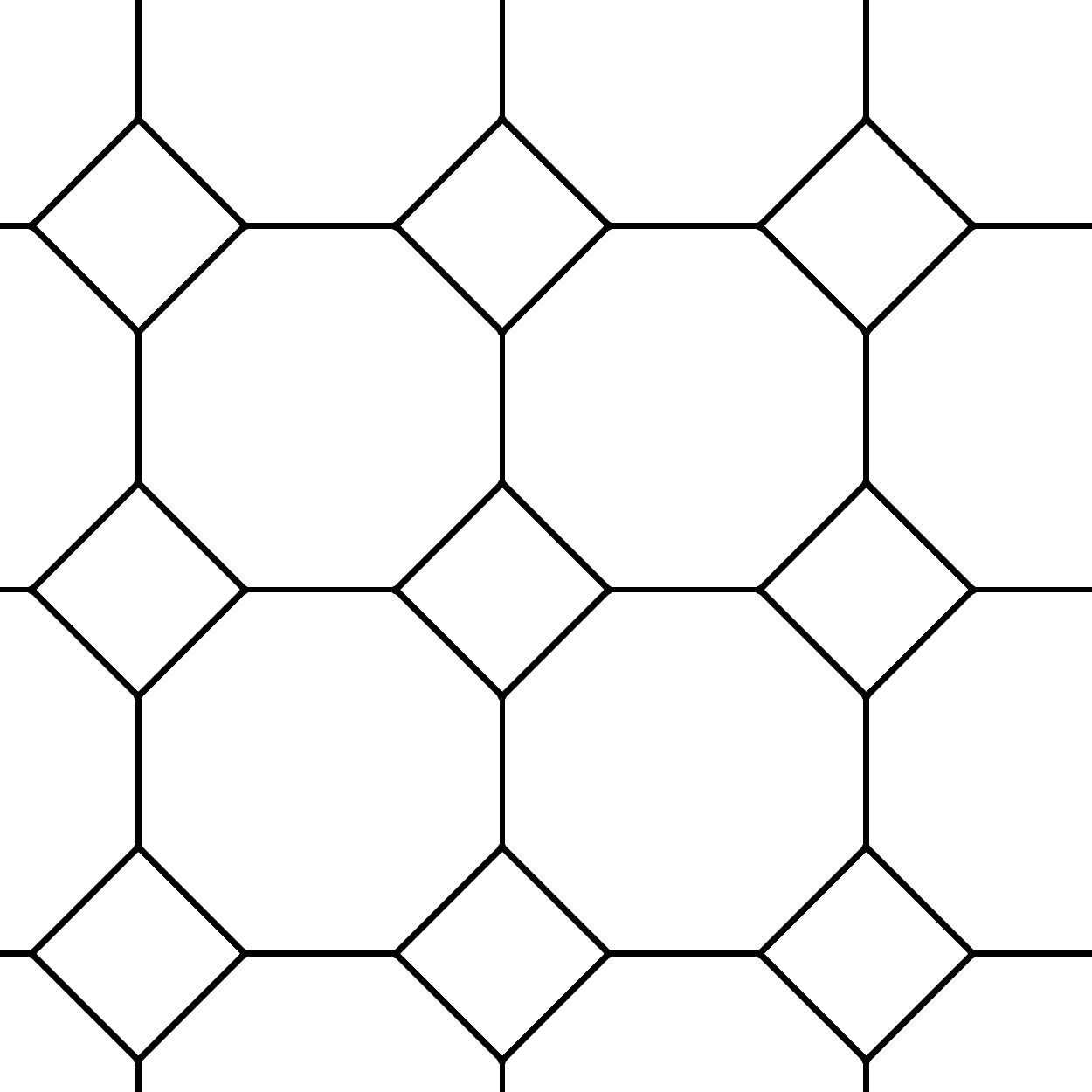} &
   \!\!$\displaystyle\frac{1}{24}{-}\frac{\arcsec(3)}{12 \sqrt{2} \pi }{+}\frac{\arcsec(3)^2}{8 \pi ^2}$\!\! &
  \multirow{2}{*}{\small $8.825563\dots$} &
  \multirow{2}{*}{\small $2.694674\dots$} &
  { \!\!$\displaystyle\frac{3}{8}{-}\frac{\arcsec(3)}{12 \sqrt{2} \pi }{+}\frac{\arcsec(3)^2}{8 \pi ^2}$\!\!} &
 { \!\!$\displaystyle\frac{9}{8}{-}\frac{\arcsec(3)}{4 \sqrt{2} \pi }{+}\frac{3 \arcsec(3)^2}{8 \pi ^2}$\!\!} &
  { \!$\displaystyle\frac{25}{16}{-}\frac{\arcsec(3)}{8 \sqrt{2} \pi }{+}\frac{3 \arcsec(3)^2}{16 \pi ^2}$\!} \\[12pt]
 & &
  \small $0.037769\dots$ &
  & &
  \small $0.371102\dots$ &
  \small $1.113307\dots$ &
  \small $1.556654\dots$ \\
\hline
&&&&&&&\\[\rowgap]
  \multirow{2}{*}{tetrakis square} &
  \lattice{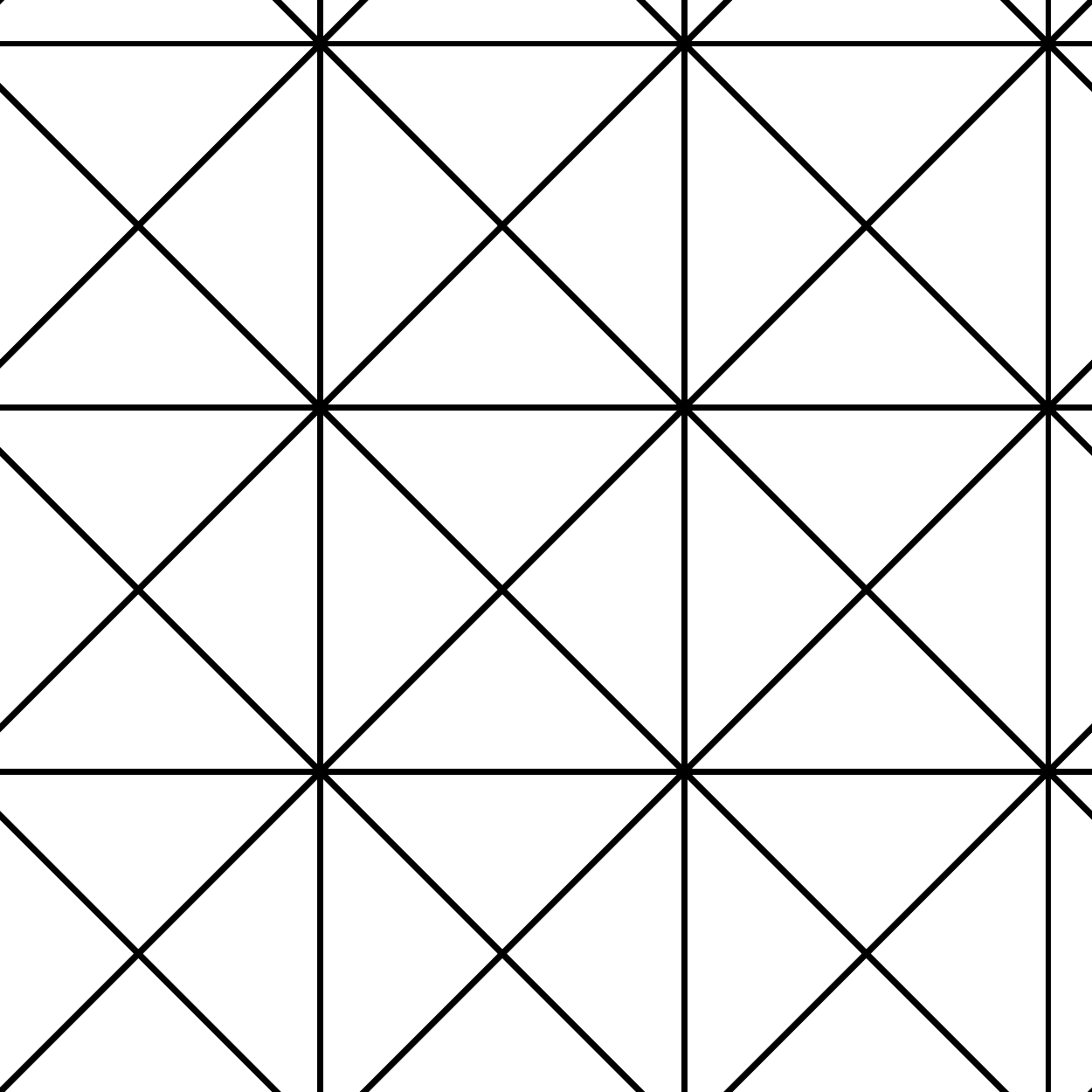} &
   \!$\displaystyle\frac{1}{8}{-}\frac{\arcsec(3)}{12 \sqrt{2} \pi }{+}\frac{\arcsec(3)^2}{16 \pi ^2}$\! &
  \multirow{2}{*}{\small $5.978703\dots$} &
  \multirow{2}{*}{\small $3.594878\dots$} &
  {\!\!$\displaystyle\frac{7}{24}{-}\frac{\arcsec(3)}{12 \sqrt{2} \pi }{+}\frac{\arcsec(3)^2}{16 \pi ^2}$\!\!} &
 { \!\!$\displaystyle\frac{7}{4}{-}\frac{\arcsec(3)}{2 \sqrt{2} \pi }{+}\frac{3 \arcsec(3)^2}{8 \pi ^2}$\!\!} &
  { \!$\displaystyle\frac{27}{8}{-}\frac{\arcsec(3)}{4 \sqrt{2} \pi }{+}\frac{3 \arcsec(3)^2}{16 \pi ^2}$\!} \\[12pt]
  & &
  \small $0.111507\dots$ &
  & &
  \small $0.278174\dots$ &
  \small $1.669041\dots$ &
  \small $3.334521\dots$ \\[2pt]
\hline
\end{tabular}
}
\vspace*{6pt}
\end{center}
\caption{Unicycle, loop-erased random walk, uniform spanning tree, and sandpile parameter values for different lattices.}
\label{values}
\end{sidewaystable}

We give a new method for computing sandpile densities of planar graphs
which is simpler and readily applies to other planar lattices.  The
calculations are elementary and require only modest background on
sandpiles and spanning trees.
The main ingredients are a
combinatorial use of planar duality and an explicit counting formula.
For the square lattice, all the intermediate expressions are rational,
and essentially depend only on simple symmetry arguments.  For the
triangular and honeycomb lattices, the sandpile densities were
computed by Kenyon and Wilson and determined to be $10/3$ and $37/24$
respectively \cite{KW4}, though these computations involved
intermediate values containing $\sqrt{3}/\pi$.  Using our new method
together with the symmetries of these lattices, we can easily recover
the $10/3$ and $37/24$ values.  There are other lattices with a high
degree of symmetry, such as the kagom\'e lattice, the dice lattice, or
the Fisher lattice, for which one can see without doing any
calculations that the sandpile density must be a rational number, and
it is not much work to compute those numbers ($13/6$, $17/8$, and
$959/600$).  For the square-octagon lattice, the sandpile density is
transcendental, but can be expressed in terms of an inverse
trigonometric function.  For other $\Z^2$-periodic graphs more
generally, the sandpile density is expressible in terms of simple
electrical quantities.

The sandpile density is closely related to certain
quantities in random spanning trees and related structures,
 including
the steady-state rate at which
loop-erased random walk (LERW) produces and then erases loops (the ``looping rate''),
the probability that the spanning tree path from a random vertex to infinity
passes through a neighboring vertex,
and the expected length of the cycle in a uniformly
random spanning unicycle (a connected spanning subgraph containing
exactly one cycle) \cite{PP:54,LP:54}.  Table~1 summarizes these
values for the various lattices mentioned above.

In Section~\ref{2SF}, building on earlier
work~\cite{liu-chow,myrvold}, we show how to compute the number of
two-component spanning forests in terms of electric current across
edges.  When the underlying graph is planar,
two-component spanning forests are related by duality to spanning
unicycles, which is what allows us to carry out the above
calculations.  For most of the above-mentioned lattices, the electric
current across edges can be evaluated by simple symmetry arguments.
In Sections~\ref{sec:tree-LERW}--\ref{sec:sandpile} we provide
further background explaining how spanning trees, electric networks,
loop-erased random walk, spanning unicycles, and sandpiles are all
related.  These different relations are valid for any finite graph,
but the exact computations we carry out rely on planarity.  In
Section~\ref{periodic} we discuss the infinite lattice limit and carry
out the concrete calculations.  We conclude in
Section~\ref{openproblems} with some open questions.

\section{The Matrix-Tree Theorem and spanning forests}\label{2SF}

An important tool in spanning tree calculations is the
\textit{Matrix-Tree Theorem}, which is essentially due to Kirchhoff.
We describe this theorem as follows.
Let $\G$ be a finite connected graph endowed with a weight
function~$w$, giving to any oriented pair of vertices $(u,v)$ a weight
$w_{u,v}$ with the convention that $w_{u,v}=0$ if $uv$ is not an edge.
When $w_{u,v} = w_{v,u}$, as first considered by Kirchhoff,
this may be viewed as an electrical network with conductance $w_{u,v}$
on the resistor $uv$.  The graph Laplacian $\Delta$ of~$\G$ is the
matrix defined by $\Delta_{u,v}=-w_{u,v}$ when $u\neq v$, and
$\Delta_{u,u}=\sum_{v\neq u} w_{u,v}$.

\begin{figure}[t!]
\[
\hspace{-30pt}\det
\underbrace{\begin{bmatrix}
w_{1,2}{+}w_{1,3}&-w_{1,2}&\smash{\makebox[0pt][l]{{\kern-10pt\textcolor{lightgray}{\rule[-28pt]{53pt}{37pt}}}}}-w_{1,3}\\
-w_{2,1}&w_{2,1}{+}w_{2,3}&-w_{2,3}\\
\smash{\makebox[0pt][l]{{\kern-8pt\textcolor{lightgray}{\rule[-4pt]{140pt}{13pt}}}}}
-w_{3,1}&-w_{3,2}& w_{3,1}{+}w_{3,2}
\end{bmatrix}}_{\text{graph Laplacian $\Delta$}}
= w_{1,2}w_{2,3} + w_{2,1}w_{1,3} + w_{1,3} w_{2,3}
\]

\definecolor {processblue}{cmyk}{0.96,0,0,0}
\newcommand{\pAB}{\path (A) edge [bend left = 25,color=processblue,thick] node[color=black,above=-8pt] {\contour{white}{$w_{1,2}$}} (B);}
\newcommand{\pAC}{\path (A) edge [bend right=-15,color=processblue,thick] node[below =-8pt,color=black,sloped] {\contour{white}{$w_{1,3}$}} (C);}
\newcommand{\pBA}{\path (B) edge [bend left = 15,color=processblue,thick] node[below=-8pt,color=black,sloped] {\contour{white}{$w_{2,1}$}} (A);}
\newcommand{\pBC}{\path (B) edge [bend right=-25,color=processblue,thick] node[below =-8pt,color=black,sloped] {\contour{white}{$w_{2,3}$}} (C);}
\newcommand{\pCA}{\path (C) edge [bend left = 25,color=processblue,thick] node[below=-8pt =0.15 cm,color=black,sloped] {\contour{white}{$w_{3,1}$}} (A);}
\newcommand{\pCB}{\path (C) edge [bend left = 15,color=processblue,thick] node[below=-9pt,color=black,sloped] {\contour{white}{$w_{3,2}$}} (B);}
\[
\hspace{10pt}
\raisebox{-0.5\height}{
\begin {tikzpicture}[-latex,auto,node distance =1.73205 cm and 1cm ,on grid , semithick , scale=0.9,
every path/.style = {transform shape},
every node/.style ={transform shape},
state/.style ={transform shape, circle, top color=white, bottom color=processblue!20, draw,processblue, text=blue , minimum width =10pt}]
\node[state] (C) {$3$};
\node[state] (A) [above left=of C] {$1$};
\node[state] (B) [above right =of C] {$2$};
\pAB \pAC
\pBC \pBA
\pCA \pCB
\end{tikzpicture}
}
\implies
\raisebox{-0.5\height}{
\begin {tikzpicture}[-latex,auto,node distance =1.73205 cm and 1cm ,on grid , semithick , scale=0.9,
every path/.style = {transform shape},
every node/.style ={transform shape},
state/.style ={transform shape, circle, top color=white, bottom color=processblue!20, draw,processblue, text=blue , minimum width =10pt}]
\node[state] (C) [color=brown] {$3$};
\node[state] (A) [above left=of C] {$1$};
\node[state] (B) [above right =of C] {$2$};
\pAB
\pBC
\end{tikzpicture}
}
+
\raisebox{-0.5\height}{
\begin {tikzpicture}[-latex,auto,node distance =1.73205 cm and 1cm ,on grid , semithick , scale=0.9,
every path/.style = {transform shape},
every node/.style ={transform shape},
state/.style ={transform shape, circle, top color=white, bottom color=processblue!20, draw,processblue, text=blue , minimum width =10pt}]
\node[state] (C) [color=brown]{$3$};
\node[state] (A) [above left=of C] {$1$};
\node[state] (B) [above right =of C] {$2$};
\pAC
\pBA
\end{tikzpicture}
}
+
\raisebox{-0.5\height}{
\begin {tikzpicture}[-latex,auto,node distance =1.73205 cm and 1cm ,on grid , semithick , scale=0.9,
every path/.style = {transform shape},
every node/.style ={transform shape},
state/.style ={transform shape, circle, top color=white, bottom color=processblue!20, draw,processblue, text=blue , minimum width =10pt}]
\node[state] (C) [color=brown]{$3$};
\node[state] (A) [above left=of C] {$1$};
\node[state] (B) [above right =of C] {$2$};
\pAC
\pBC
\end{tikzpicture}
}
\]
\caption{The Matrix-Tree Theorem for directed graphs.
When the graph Laplacian $\Delta$ has row $s$ and column $s$
removed, the determinant of the resulting matrix $\Delta_{\widehat{s}}^{\widehat{s}}$
gives the weighted sum of arborescences rooted at $s$.
See \cite{Temperley} or \cite{zeilberger} for an elegant combinatorial proof.
}\label{fig:MTT}
\end{figure}
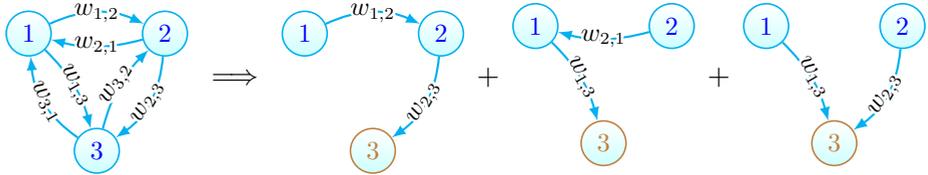

The theorem gives the weighted sum of spanning trees of $\G$ (the weight
of a tree is the product of its edge weights)
as the determinant of a submatrix of the Laplacian
$\Delta$.  Specifically, for any vertex $s$ of~$\G$, if we remove the
row and column associated with~$s$, then the determinant of the
resulting matrix gives the weighted sum of spanning trees, that is
\[F_1(\G)=\det \Delta_{\widehat{s}}^{\widehat{s}}\,,\]
where $F_1(\G)$ denotes the weighted sum of spanning trees of $\G$,
and $\Delta_{\widehat{r_1,\dots,r_k}}^{\widehat{c_1,\dots,c_k}}$ denotes the submatrix of $\Delta$
obtained by deleting rows $r_1,\dots,r_k$ and columns $c_1,\dots,c_k$.

Tutte generalized the theorem to directed graphs, and the directed
version is illustrated in Figure~\ref{fig:MTT}.  The Matrix-Tree
Theorem has been further generalized in a variety of ways
\cite{chaiken,kalai,duval-klivans-martin,Forman,Kenyon.bundle}.
Of interest to us here is a formula for counting (unrooted) spanning
forests.  Let $F_k(\G)$ denote the weighted sum of $k$-component
spanning forests of an undirected graph $\G$
(where the weight of a forest is the product of its edge weights).
Liu and Chow
\cite{liu-chow} gave a nice formula for $F_k(\G)$; the original proof
was complicated, but a short and elegant proof was given by Myrvold
\cite{myrvold}.  We shall use this formula for two-component spanning
forests, so we state and prove the formula for the special case $k=2$;
the formula for general~$k$ (given later in \eqref{eq:F_k}) and its proof are not significantly harder.
For any vertex $s$ of $\G$,
\begin{equation}\label{forestformula}
F_2(\G) = \sum_{v\neq s} \det\Delta_{\widehat{v,s}}^{\widehat{v,s}} - \sum_{\substack{u\sim v\\u,v\neq s}} w_{u,v} \det\Delta_{\widehat{u,v,s}}^{\widehat{u,v,s}} \,,
\end{equation}
where $\sum_{u\sim v}$ denotes a sum over undirected edges.

Here is the proof of \eqref{forestformula}.
For any nonempty set of vertices $S$ of $\G$,
the principal minor $\det\Delta_{\widehat{S}}^{\widehat{S}}$ gives the weighted sum of
spanning trees of the graph obtained from~$\G$ by gluing together the
vertices in $S$, or equivalently, the weighted sum of spanning forests
with $|S|$ trees in which each vertex of $S$ is in a separate tree.
The first term $\sum_{v\neq s} \det\Delta_{\widehat{v,s}}^{\widehat{v,s}}$ in \eqref{forestformula}
gives a weighted sum of two-component spanning forests in which the
tree not containing~$s$ has an extra weight which is the number of its
vertices.  The second term is a sum over three-component spanning
forests in which $u$, $v$, and $s$ are in separate trees, times the
weight of edge $(u,v)$.  But this is just a sum over two-component
spanning forests in which the tree not containing $s$ has an extra
weight which is the number of its edges.  Since any tree has one more
vertex than edge, the difference between these terms is just the
weighted sum of two-component spanning forests.

The Green's function $G$ of the graph~$\G$ with Dirichlet
boundary conditions at vertex~$s$ is given by the inverse of the
Laplacian with row~$s$ and column~$s$ removed:
\[G^{(s)}_{u,v}=\begin{cases} \big[(\Delta_{\widehat{s}}^{\widehat{s}})^{-1}\big]_{u,v} & u,v\neq s, \\ 0 & \text{$u=s$ or $v=s$.}\end{cases}\]
Since the Laplacian is symmetric, $G^{(s)}_{u,v}=G^{(s)}_{v,u}$.
The Green's function has the following electrical interpretation.  Suppose that
each edge of the graph is a conductor with conductance given by its
weight.  When one unit of current is inserted at $u$ and extracted at
$s$, it gives the \textit{voltage drop\/} between the vertices $v$ and~$s$.
Usually the vertex~$s$ is
suppressed from the notation.

As discussed in \cite{myrvold}, the forest
formula~\eqref{forestformula} can also be expressed, using Jacobi's
identity, in terms of the Green's function as
\begin{align*}
 \frac{F_2(\G)}{F_1(\G)}
 &= \sum_{v\neq s} \det G_{v}^{v} - \sum_{\substack{u\sim v\\u,v\neq s}} w_{u,v} \det G_{u,v}^{u,v}  \\
 &= \sum_{v\neq s} G_{v,v} - \sum_{\substack{u\sim v\\u,v\neq s}} w_{u,v} \big[G_{u,u} G_{v,v} - G_{u,v}^2\big]\,,
\end{align*}
where $G_{r_1,\dots,r_k}^{c_1,\dots,c_k}$ denotes the submatrix of $G$
consisting of rows $r_1,\dots,r_k$ and columns $c_1,\dots,c_k$.

Our aim is to do calculations for infinite lattices, or equivalently,
for large subgraphs in the limit where the subgraphs tend to the
infinite lattice.  In the above formula there is cancellation --- there
are quantities being added and subtracted --- and this cancellation
becomes more significant as the graphs become large, since the Green's
function diverges.  To take a limit as the graphs tend to the infinite
lattice, we re-express this formula as a sum of positive terms.

It is convenient to work with the Green's function difference
\[ A^{(s)}_{u,v} = G^{(s)}_{u,u}-G^{(s)}_{u,v}\,.\]
$A^{(s)}_{u,v}$ gives the voltage drop between vertices $u$ and $v$
when one unit of current
is run through the network from~$u$ to~$s$, so in fact $A^{(s)}_{u,v} = G^{(u)}_{s,v}$.
If $(u,v)$ is an edge, since $w_{u,v}$ is
its conductance, by Ohm's law the electric current flowing across the edge is $w_{u,v} A^{(s)}_{u,v}$.
 From this electrical interpretation, it follows that
\[ \sum_v w_{u,v} A_{u,v} = \begin{cases} 1 & u\neq s, \\ 0 & u=s.\end{cases}\]
Using this way of writing $1$ or $0$, we rewrite the for formula for $F_2(\G)/F_1(\G)$ as
\begin{align}
 \frac{F_2(\G)}{F_1(\G)}
 &= \sum_{u\sim v} w_{u,v} [G_{u,u} A_{u,v} + G_{v,v} A_{v,u} - G_{u,u} G_{v,v} + G_{u,v}^2]\,, \notag\\
\intertext{where we no longer exclude edges incident to $s$ from the sum, since those terms
contribute $0$ anyway.  This formula may be further rewritten as}
\label{two-tree}
 \frac{F_2(\G)}{F_1(\G)} &= \sum_{u\sim v} w_{u,v} \bigg[\big(A^{(s)}_{u,v} - A^{(s)}_{v,u}\big)^2 + A^{(s)}_{u,v} A^{(s)}_{v,u} \bigg] \,.
\end{align}
Equation~\eqref{two-tree} holds for any finite undirected weighted
graph $\G$ and vertex $s$ of $\G$.  It is ideal for our purposes.  All
the terms are positive, so there is no cancellation, and in many cases
of interest it is easy to evaluate
the $A_{u,v}$'s for neighboring vertices.

\section{Cycle-rooted spanning trees and loop-erased random walk}
\label{sec:tree-LERW}

There is a natural Markov chain on spanning trees which is as
important to their analysis as the Matrix-Tree Theorem.  For a finite
weighted directed graph $\G$, an \textit{arborescence\/} is a spanning
tree of $\G$ in which all edges are directed towards some root vertex.
If we adjoin an outgoing edge from the root, the result is an
oriented cycle to which every vertex is connected via a directed path,
and is called an oriented \textit{cycle-rooted spanning tree\/}
(CRST).  It is useful to place a mark at the root of the arborescence
before adjoining the extra edge, so that the oriented CRST has a
marked vertex on its cycle.

\begin{figure}[b!]
\centering
\includegraphics[width=\textwidth]{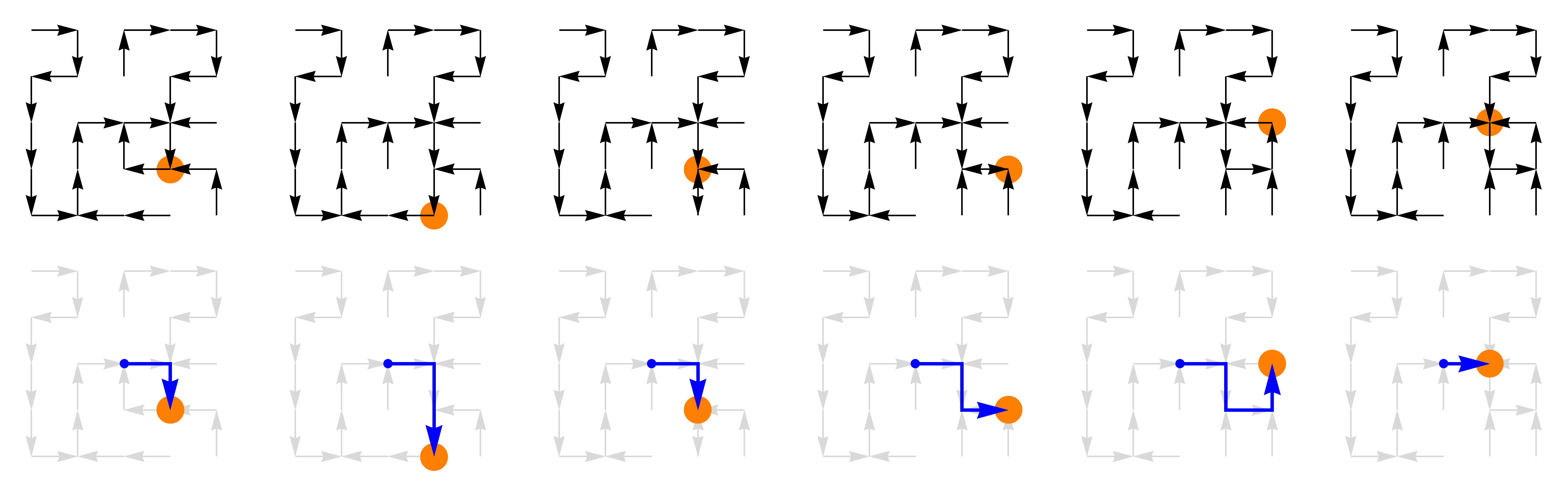}
\caption{A few steps of the Markov chain on marked oriented cycle-rooted spanning trees on
a $5\times 5$ square grid, shown in the top row.  The bottom row shows the path within the CRST from a fixed starting
vertex leading to the mark; this path evolves according to a loop-erased random walk.
The LERW erases a loop precisely when
the LERW intersects the CRST cycle only at the mark; the erased loop is the CRST cycle.
}
\label{MarkovChain}
\end{figure}

In the canonical probability distribution on marked oriented CRST's,
each configuration occurs with a probability that is proportional to the
product of the weights of its edges.  Consider the following Markov
chain on marked oriented CRST's (see Figure~\ref{MarkovChain}).  At each step, the
Markov chain erases the outgoing edge from the marked vertex, picks a
new random outgoing edge from that vertex with probability
proportional to the edge weights, and then moves the mark forward one
step along the new cycle.  Picking a new outgoing edge from the marked
vertex preserves the canonical probability distribution, as does
sliding the mark forward, so the canonical probability distribution is
an invariant distribution of the Markov chain.

We now show that the canonical distribution is the unique invariant
distribution provided the graph is strongly connected, meaning that there is a
directed path from any vertex to any other vertex.  To avoid periodicity issues, we
consider first the Markov chain run in continuous time at rate~$1$.
We run two copies of the Markov
chain starting from any two marked oriented CRST's, where the two
copies of the Markov chain are independent of one another unless
by chance 
their marks are on the same vertex,
in which case the mark
follows the same trajectory in both Markov chains.  
Since the graph is strongly connected, with
probability~$1$, the marks in the two copies of the Markov chain
will eventually be at the same vertex.
Almost surely the mark subsequently visits and departs each vertex of the
graph, and after this time both Markov chains will be in the same marked oriented CRST.
From this it is easy to see that the Markov chain can have only one
invariant distribution: otherwise we could start the two copies of the
Markov chain from random samples from the two invariant distributions,
both of which are preserved by the Markov chain, and yet with
probability~$1$ the two configurations become equal.  Thus the
canonical probability distribution on marked oriented CRST's is the
unique stationary distibution of the Markov chain.  Since any stationary
distribution for the discrete time Markov chain is also invariant for the
continuous time rate-$1$ Markov chain, the canonical distribution is also the
unique stationary distribution of the discrete time Markov chain.

If we erase the outgoing edge from the mark in the marked oriented
CRST, the result is a random arborescence with probability
proportional to the product of its edge weights times the weighted
outdegree of the root.  If we instead run the Markov chain so that
the rate at which the mark moves is the weighted degree of the vertex
where it is located, then the stationary distribution of the
arborescence becomes just the product of its edge weights.  If we
consider just the location of the mark, it does a random walk on the
underlying graph, with transition rates given by the edge weights,
i.e., it moves according to the continuous-time Markov chain defined by the weighted
graph.  This implies the Markov Chain Tree Theorem, which gives a
Markov chain's stationary distribution as being proportional to the
weighted sum of arborescences rooted at different vertices, and which
Aldous has called ``perhaps the most frequently rediscovered result in
probability''.

We can also consider the path from a fixed starting vertex $X_0$ to
the mark within the marked oriented CRST.  This path evolves according
to a loop-erased random walk (LERW),
which is a process that was introduced and first studied by Lawler
\cite{lawler-limic}.  \textit{Loop-erased random walk\/} is obtained
from random walk by erasing loops as they appear.  If
$X_0,X_1,\dots,X_t$ are the vertices of a random walk, then let $t'$ be
the largest index for which $X_{t'}=X_0$.  The loop erasure of $X_0,\dots
X_t$ is $X_0$ followed by the loop erasure of $X_{t'+1},\dots X_t$.
Consequently, the path within a uniform spanning tree from a vertex
to the root is a loop-erased random walk; this fact was 
first noted and used by Pemantle \cite{Pemantle}.
Wilson~\cite{Wilson} described further connections between random spanning
trees and random walk, giving an exact sampling algorithm
(this is how Figure~\ref{fig:UST} was produced) with implications for the
analysis of spanning trees \cite{BLPS,Schramm:SLE,LSW,LP:book}.

\section{Looping rate}
\label{sec:looping-rate}

Comparing the loop-erased random walk to the Markov chain on marked oriented cycle-rooted spanning trees, we see that the LERW creates and erases a loop precisely when the LERW first reaches the oriented CRST's cycle at the mark (see Figure~\ref{MarkovChain}).  Thus the steady state rate $\rho$ at which the (discrete time) LERW creates and erases loops is
\[
\rho=\text{discrete-time LERW looping rate}=\frac{\text{weighted sum of oriented CRST's}}{\text{weighted sum of marked oriented CRST's}}\,.
\]
The letter $\rho$ is mnemonic, since
it resembles a path with a loop at the end.  Each erased loop has size $2$ or
more.  We also let $\tau$ denote the steady-state rate at which loops
of size at least $3$ are produced.
By similar reasoning,
\[
\tau=\frac{\text{weighted sum of oriented CRST's with cycle length $\geq 3$}}{\text{weighted sum of marked oriented CRST's}}\,.
\]
Since the above formulas do not depend on the start vertex $X_0$,
it follows \textit{a posteriori\/} that these looping rates $\rho$
and $\tau$ are independent of the start vertex.

As we shall see, the sandpile
density is related to $\rho$, the difference between~$\rho$ and~$\tau$
is easy to compute, and $\tau$ is closely related to spanning
unicycles, which for planar graphs are equivalent to two-component
spanning forests on the dual graph.

Suppose that the underlying graph~$\G$ is undirected.  Given a marked oriented CRST, we can separate the mark and its outgoing edge from the CRST.  The result is a spanning tree~$T$ together with an edge~$e$ and a distinguished endpoint of $e$.  In the reverse direction, a spanning tree~$T$ and an edge~$e$ with distinguished endpoint can be combined to form a marked oriented cycle-rooted spanning tree, where the mark is at the distinguished vertex, the cycle is oriented the direction of $e$ from the mark, and the other edges are oriented towards the cycle.
Thus for undirected graphs,
\begin{equation} \label{eq:rho-tau}
\begin{split}
\rho-\tau
&=\frac{\text{weighted sum of oriented CRST's with cycle length 2}}{\text{weighted sum of trees} \times 2\times \text{weighted sum of edges}} \\
&=\frac12\Pr[\text{random edge }e\in \text{random tree }T]\,,
\end{split}
\end{equation}
where the random edge $e$ and random tree $T$ are chosen according to the edge weights.

For unweighted undirected graphs this probability is trivial to evaluate, since
regardless of what the tree~$T$ is, the conditional probability that a random edge
is in $T$ is $(|V|-1)/|E|$, where $V$ and $E$ are the vertex and edge sets of the graph $\G$, so
\begin{equation} \label{eq:r-e-in-T-undirected}
  \Pr[\text{random edge }e\in \text{random tree }T] = \frac{|V|-1}{|E|}
\quad\quad\text{(unweighted graphs)}\,.
\end{equation}

For weighted undirected graphs, we can compute this probability using
the connections between spanning trees and random walk,
and between random walk and electric networks.
Let $u\leadsto x\to y\leadsto s$ denote the event that the path from $u$
to $s$ in a random spanning tree includes the edge $x\sim y$ in the direction
from $x$ to $y$.  The path from $u$ to $s$ is a loop-erased random walk.
For undirected graphs, each erased loop is equally likely to have been traversed in
either direction.  Thus $\Pr[u\leadsto x\to y\leadsto s] - \Pr[u\leadsto y\to x\leadsto s]$
is the expected algebraic number of traversals of the edge $x\sim y$, i.e.,
traversals of $(x,y)$ minus traversals of $(y,x)$, made by a random walk started from
$u$ and stopped at $s$.  The Green's function $G^{(s)}_{u,x}$ gives the expected time
spent at $x$ by a continuous time random walk started at $u$ and stopped at $s$.
Thus the expected algebraic number of traversals of edge $x\sim y$
is $w_{u,v} (G^{(s)}_{u,x} - G^{(s)}_{u,y})$, which has the interpretation of the
current flowing across edge $(x,y)$ when one unit of current is inserted at $u$
and extracted at $s$.  (See \cite{doyle-snell} for further background on random walks
and electric networks.)
Hence
\begin{equation} \label{current}
\Pr[u\leadsto x\to y\leadsto s] - \Pr[u\leadsto y\to x\leadsto s]  = w_{u,v} (G^{(s)}_{u,x} - G^{(s)}_{u,y})\,.
\end{equation}
Taking $x=u$ and $y=v$ gives
$\Pr[u\to v\leadsto s] = w_{u,v} A^{(s)}_{u,v}$, from which it follows that
\begin{equation} \label{uv-in-T}
\Pr[(u,v)\in T]=w_{u,v}(A_{u,v}+A_{v,u})\,,
\end{equation}
and that when the edge $e$ is chosen at random according to the weights $w$,
\begin{equation} \label{eq:r-e-in-T}
  \Pr[e\in T]
 = \frac{\displaystyle\sum_{u\sim v} w_{u,v}^2 (A_{u,v} + A_{v,u})}{\displaystyle\sum_{u\sim v} w_{u,v}}
  \quad\quad\quad\text{(undirected graphs)}.
\end{equation}

No good formula or efficient algorithm is known for counting spanning
unicycles of a general graph.  But for \textit{planar\/} graphs, the
dual of a spanning unicycle is a two-component spanning forest, for
which \eqref{two-tree} gives a weighted count.
If $\G$ is embedded in the plane, let $\G^*$ denote its planar dual
($\G^*$ depends on the embedding).  Each edge $e$ of~$\G$ has a dual edge $e^*$ with weight
$w(e^*)=1/w(e)$.
Observe that by planar duality
\[\text{weighted sum of unicycles of $\G$}=F_2(\G^*)\prod_{e\in E}w(e) \,,\]
and
$F_1(\G)=F_1(\G^*)\prod_{e\in E}w(e)$ .  Using \eqref{two-tree}
we therefore obtain
\begin{equation}\label{tau}
\tau
= \frac{\displaystyle\sum_{u^*\sim v^*} w_{u^*,v^*}\left(A^{(s^*)}_{u^*,v^*}A^{(s^*)}_{v^*,u^*}+\left(A^{(s^*)}_{u^*,v^*}-A^{(s^*)}_{v^*,u^*}\right)^2\right)}{\displaystyle\sum_{u^*\sim v^*} 1/w_{u^*,v^*}}\,,
\end{equation}
where the sums are over adjacent faces $u^*$ and $v^*$ of $\G$, i.e., adjacent vertices in $\G^*$.

\section{Mean loop length and neighboring spanning tree ancestors}
For finite graphs, the expected LERW loop length is $1/\rho$,
since after~$N$ steps there are about $\rho N$ loops which altogether
contain about $N$ edges.

We define $\lambda$ to be the \textit{mean unicycle loop length}, that is,
the expected number of edges in the cycle of a $w$-random spanning
unicycle.  In terms of LERW, we see that after $N$ steps there are
about $\tau N$ long loops, which altogether contain about
$N-2(\rho-\tau)N$ edges, so
\begin{equation}\label{lambda}
 \quad\lambda= \frac{1-2(\rho-\tau)}{\tau} = \frac{\Pr[e\notin T]}{\tau} \quad\quad\quad\text{(undirected graphs)},
\end{equation}
where $e$ is a $w$-random edge and $T$ is a $w$-random spanning tree.
For unweighted graphs this relation appears in \cite{LP:54}.

There is another interpretation for the looping rate which is discussed in \cite{LP:54}.
Recall that for undirected graphs a random marked oriented CRST is the union of a random spanning tree $T$ with an independent random directed edge~$e$.  Since $\rho$ is the probability that the marked oriented CRST has its mark on the path from $s$ to the cycle,
\begin{align*}
\rho & =\sum_{u}\sum_{v\sim u} \Pr[e=(u,v)] \times \Pr[s\leadsto v\leadsto u\text{ in $T$}]\,.
\end{align*}
Thus the (weighted) mean number of neighboring ancestors of a uniformly
chosen random vertex $u$ in a random spanning arborescence $T$ rooted at $s$ is
\[
\E\, \sum_{v\sim u} w_{u,v} 1_{\{u\leadsto v\leadsto s\text{ in $T$}\}} = \delta\rho\,,
\]
where we let $\delta=2\sum_{u\sim v}w_{u,v}/|V|$ denote the mean weighted degree.

\section{Sandpile density}
\label{sec:sandpile}

We outline here the key facts we use relating sandpiles and spanning
trees, which allow us to use what we know about LERW to compute the
density of sand in recurrent sandpiles.  See \cite{HLMPPW} or
\cite{jarai:sandpile} for a more in-depth introduction to sandpiles
for mathematicians.

\begin{figure}[b]
\includegraphics[width=\textwidth]{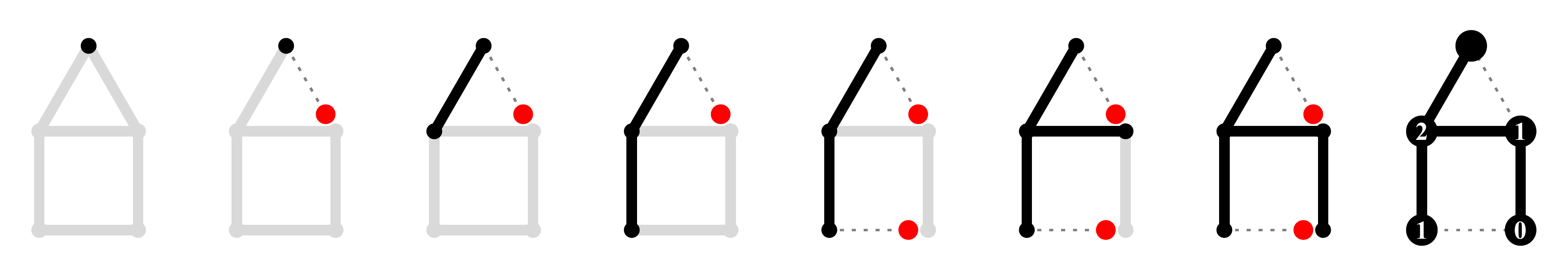}
\caption{The exploration process querying edges and converting the responses into both a spanning tree and a sandpile.
To convert a spanning tree into a sandpile, the spanning tree responds yes to edges that it contains.
To convert a sandpile to a spanning tree, the sandpile responds yes once for each grain of sand the source of the edge contains.}
\label{fig:bijection}
\end{figure}

An essential tool is the bijective map between recurrent sandpiles and spanning
trees which was first exhibited by Majumdar and Dhar
\cite{Majumdar-Dhar}.  The map between trees and sandpiles is not
canonically unique, and since their work, several variations and
generalizations have been published (notably
\cite{cori-le-borgne,bernardi,jarai-werning:measures}), with different
mappings being useful for different purposes.  The maps from trees to
sandpiles correspond to tree exploration processes.
We describe such an exploration process
as follows, see Figure~\ref{fig:bijection}.
Imagine that there
is an arborescence which is hidden, except for the root~$s$, which is the
initial ``current tree.''  For any edge leading from a vertex $u$ not
in the current tree to a vertex $v$ in the current tree, we may query
if that edge $(u,v)$ is in the tree; if so, then vertex~$u$ and edge
$(u,v)$ get adjoined to the current tree, and otherwise~$u$ gains a
mark.  (The main difference between the different variations of the map
 is the
choice of which such edge to query next.)  The final amount of sand at
a vertex is its out-degree minus one minus the number of its marks.
The resulting sandpile configuration is
non-negative and stable.  It is also clear that the map from trees to
sandpiles is one-to-one, since given the sandpile configuration, we
can run the same exploration process (as in
Figure~\ref{fig:bijection}), and use the sandpile heights to determine
which edge queries resulted in a yes or no answer.  It may not be
immediately clear that the resulting sandpile configuration is
recurrent, that every recurrent sandpile configuration may be
obtained from a tree in this way, and that the stationary distribution
on sandpiles is uniformly distributed on the recurrent sandpiles.
For this we refer the reader to the exposition \cite{HLMPPW}.

We make use of a formula
which, for undirected graphs, expresses the
number of sandpiles with different amounts of sand in terms of the
Tutte polynomial.  
The Tutte polynomial of an undirected graph $\G$
(with edge set $E$ and vertex set $V$) is a polynomial in two
variables defined by
\begin{equation}
T_\G(x,y)=\sum_{E'\subseteq E}(x-1)^{k(E')-k(E)}(y-1)^{k(E')+|E'|-|V|}\,,
\end{equation}
where $k(E')$ is the number of connected components of the spanning subgraph of $\G$ with edge set $E'$.

Biggs defined the \textit{level\/} of a sandpile configuration to be
the amount of sand shifted down by $|E|-\delta_s$, and showed that $0 \leq \text{level}
\leq |E| - |V| +1$.  Moreover, these bounds are tight.  Biggs conjectured
and Merino proved \cite{merino-lopez} that for connected graphs~$\G$, the generating function of recurrent
sandpiles by level is
\begin{equation}\label{merino}
\sum_{\substack{\text{recurrent}\\\text{sandpiles $\sigma$}}} y^{\level(\sigma)} = T_\G(1,y)\,.
\end{equation}
Merino proved \eqref{merino} by induction on the number of edges of
the graph.  Cori and Le~Borgne \cite{cori-le-borgne} gave a bijective
proof of \eqref{merino} using another formula for the Tutte polynomial,
which expresses $T_\G(x,y)$ as a weighted sum of spanning trees, together with
the correspondence between sandpiles and spanning trees.

The expected amount of sand can be expressed in terms of the number of
unicycles, and more generally, the $j$th binomial moment of the sandpile level can be
obtained by differentiating \eqref{merino} $j$ times and evaluating at $y=1$:
\begin{equation} \label{T'(1,y)}
\sum_{\substack{\text{recurrent}\\\text{sandpiles $\sigma$}}}
  \binom{\level(\sigma)}{j}
  = \left.\frac{1}{j!}\frac{d^j}{dy^j} T_\G(1,y)\right|_{y=1}
  = \parbox{1.6in}{\centering\# connected subgraphs of $\G$ with $|V|+j-1$ edges.}
\end{equation}
Comparing the cases $j=1$ and $j=0$, we see that
for a random recurrent sandpile, the expected level is
\[
\E[\level(\sigma)] = \frac{\text{\# unicycles of $\G$}}{\text{\# spanning trees of $\G$}} = \tau\times |E|\,,
\]
and consequently, the sandpile density $\bar\sigma = \frac{1}{|V|}\sum_v \E[\sigma(v)]$ is
\begin{equation}\label{eq:sand-density}
 \text{sandpile density} =  \bar\sigma = \frac{\delta\rho+\delta-1}{2} - \frac{\delta_s-1/2}{|V|}\,,
\end{equation}
where $\delta$ is the weighted average degree.

\section{Periodic planar lattices}
\label{periodic}

\begin{figure}[t]
\centering{
\hfill
\includegraphics[height=120pt,rotate=90]{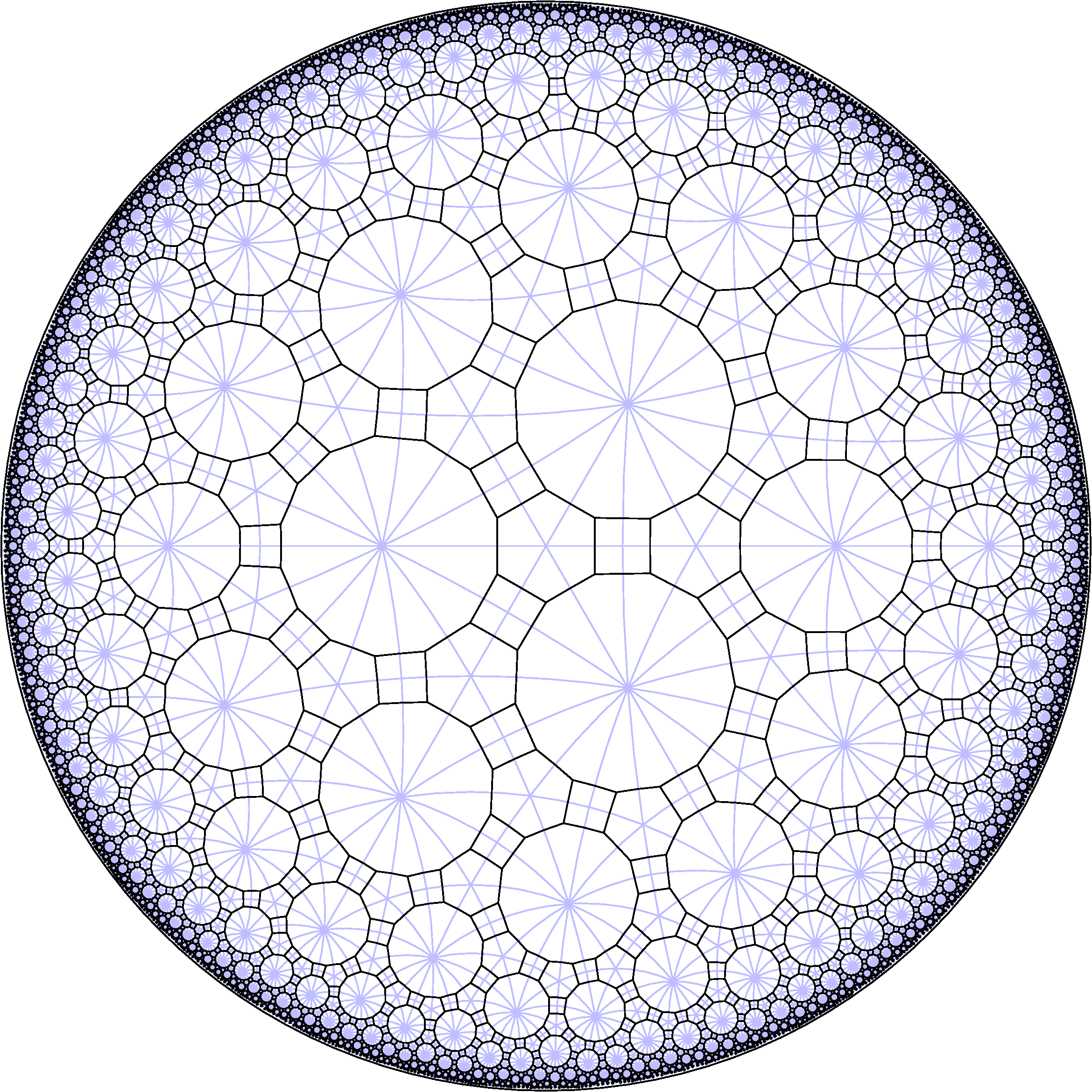}\hfill
\includegraphics[height=120pt]{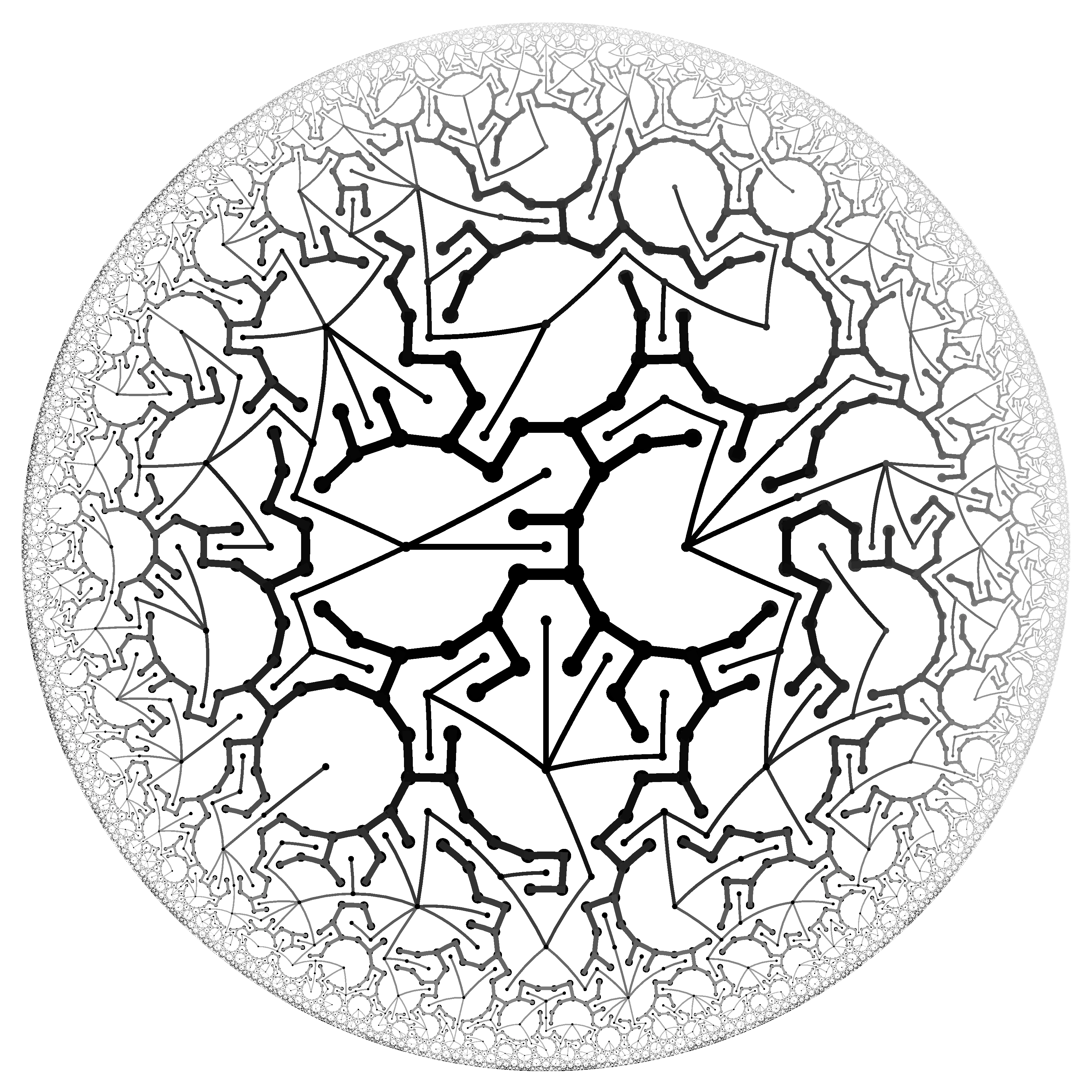}\hfill}
\caption{A hyperbolic planar lattice with its dual lattice (left),
 and a ``wired'' uniform spanning forest on one lattice together
 with its dual ``free'' uniform spanning forest on the dual lattice (right).
 (The left figure was made by Don Hatch's \texttt{HyperbolicApplet} program,
the right figure was made by Russell Lyons \cite[Figs.~6.1 and~10.3]{LP:book}.)}
\label{fig:hyperbolic}
\end{figure}

In this section, we show how to compute the looping rate and sandpile
density for periodic planar lattices.  To a large extent, the formulas
for the infinite lattices follow from the finite-graph formulas, but
some care is needed when taking the limit where the graph tends to the
infinite lattice.  For example, there are hyperbolic planar lattices,
such as the one shown in Figure~\ref{fig:hyperbolic}, for which there
is more than one uniform spanning tree measure, and more than one
infinite sandpile measure.  For finite graphs approximating the
hyperbolic lattice, in a certain sense the boundary is too big to be
negligible, and different boundary conditions lead to different
limiting uniform spanning forest measures and different sandpile
measures.  To make sense of quantities such as ``the sandpile
density,'' we would like to know that there is only one canonical
infinite sandpile measure.

Here we shall limit our attention to $\Z^2$-periodic connected planar
graphs in which the fundamental domain has a finite number of
vertices.  Each of the lattices listed in Table~1 is of this type.
Graphs of this type are recurrent (which follows e.g., from Rayleigh
monotonicity \cite{doyle-snell}).  Any recurrent graph has a unique limiting
uniform spanning forest measure, and almost surely the spanning forest
contains a single tree \cite{BLPS}.

Since the uniform spanning forest on the dual lattice is also unique
and almost surely contains a single tree, it follows that the spanning tree on the
primal lattice almost surely has one end, i.e., almost surely it does not contain
a doubly infinite path.

Whenever the uniform spanning forest almost surely contains a single tree with one end,
it is clear that there is a unique limiting sandpile measure that results from exploring
the tree.  (In cases where the spanning forest contains multiple trees, depending on the details
of the tree exploration process near the boundary, the USF trees could be explored in different
orders, resulting in different sandpiles.  However, using a carefully selected tree exploration rule,
J\'arai and Werning showed that whenever the USF almost surely contains one-ended trees,
there is a unique limiting sandpile measure \cite{jarai-werning:measures}.)

Let $\G_n$ be the graph obtained from the $\Z^2$-periodic lattice by merging
all vertices outside an $n\times n$ block of fundamental domains
(this gives ``wired boundary conditions'').
The sequence of graphs $(\G_n)_{n\geq 1}$ converges to the
lattice in the sense of Benjamini and Schramm
\cite{benjamini-schramm}, which is to say that for any distance $j>0$,
the $j$-neighborhood of a random vertex of $\G_n$ converges in
distribution as $n\to\infty$ to the $j$-neighborhood of a random
vertex in the fundamental domain.  (In contrast, the hyperbolic
lattice is not a Benjamini-Schramm limit of planar graphs
\cite{benjamini-schramm}.)  Since the USF in the lattice has one tree,
the spanning tree path connecting two random vertices in $\G_n$, when
restricted to a neighborhood of one of the vertices, converges in law
to the LERW from a random vertex in the lattice to $\infty$.  In
particular, the looping rate of~$\G_n$ converges to the looping rate
of the lattice, as does the distribution of the erased LERW loops.
Comparing the tree exploration process on $\G_n$ to that on the
lattice, since almost surely the USF on the lattice has one tree with one end,
the distribution of sand around a random vertex of~$\G_n$ converges to
the distribution of sand around a random vertex of the lattice, and
since the sand at each vertex of $\G_n$ is bounded, the rare vertices
of $\G_n$ with atypical neighborhoods can be ignored, and the density of
sand of~$\G_n$ converges to the lattice sandpile density.
Consequently,
the finite graph formulas relating the LERW looping rate and other
graph parameters also hold in the setting of $\Z^2$-periodic lattices.

Furthermore, since the USF is unique and almost surely has a single one-ended tree,
\eqref{current} implies that there is a unique limit as $n\to\infty$ for the current flowing across any
edge in the lattice.  Consequently the voltage drop $A^{(s)}_{u,v}$ has a unique limit,
which we call the \emph{potential kernel}.
Since it is unique, the potential kernel inherits all the symmetries of the lattice.

We compute $\tau$ using \eqref{tau}.
For unweighted lattices, we do not explicitly compute $\rho$, since
$\rho-\tau$ is expressible in terms of the average degree~$\delta$ via
\eqref{eq:rho-tau} and \eqref{eq:r-e-in-T-undirected}.
For weighted lattices, we compute $\rho$ using \eqref{eq:rho-tau} and \eqref{eq:r-e-in-T}.

\vspace{9pt}

At this point we can start calculating.  For the square lattice, by
symmetry the potential kernel across neighboring vertices is $1/4$, so
for its dual, also the square lattice, we have
\[ \tau\text{\tiny(square)} = \frac{1}{4} \times \frac{1}{4} + \left(\frac14-\frac14\right)^2 = \frac{1}{16}\,.\]

For the triangular lattice, the dual is the honeycomb lattice, for which by symmetry the current flowing across edges is $1/3$,
\[ \tau\text{\tiny(triangular)} = \frac{1}{3} \times \frac{1}{3} + \left(\frac13-\frac13\right)^2 = \frac{1}{9}\,.\]

For the honeycomb lattice, the dual is the triangular lattice, and we have
\[ \tau\text{\tiny(honeycomb)} = \frac{1}{6} \times \frac{1}{6} + \left(\frac16-\frac16\right)^2 = \frac{1}{36}\,.\]

\begin{figure}[h]
\centering{\hfill
\includegraphics[width=0.35\textwidth]{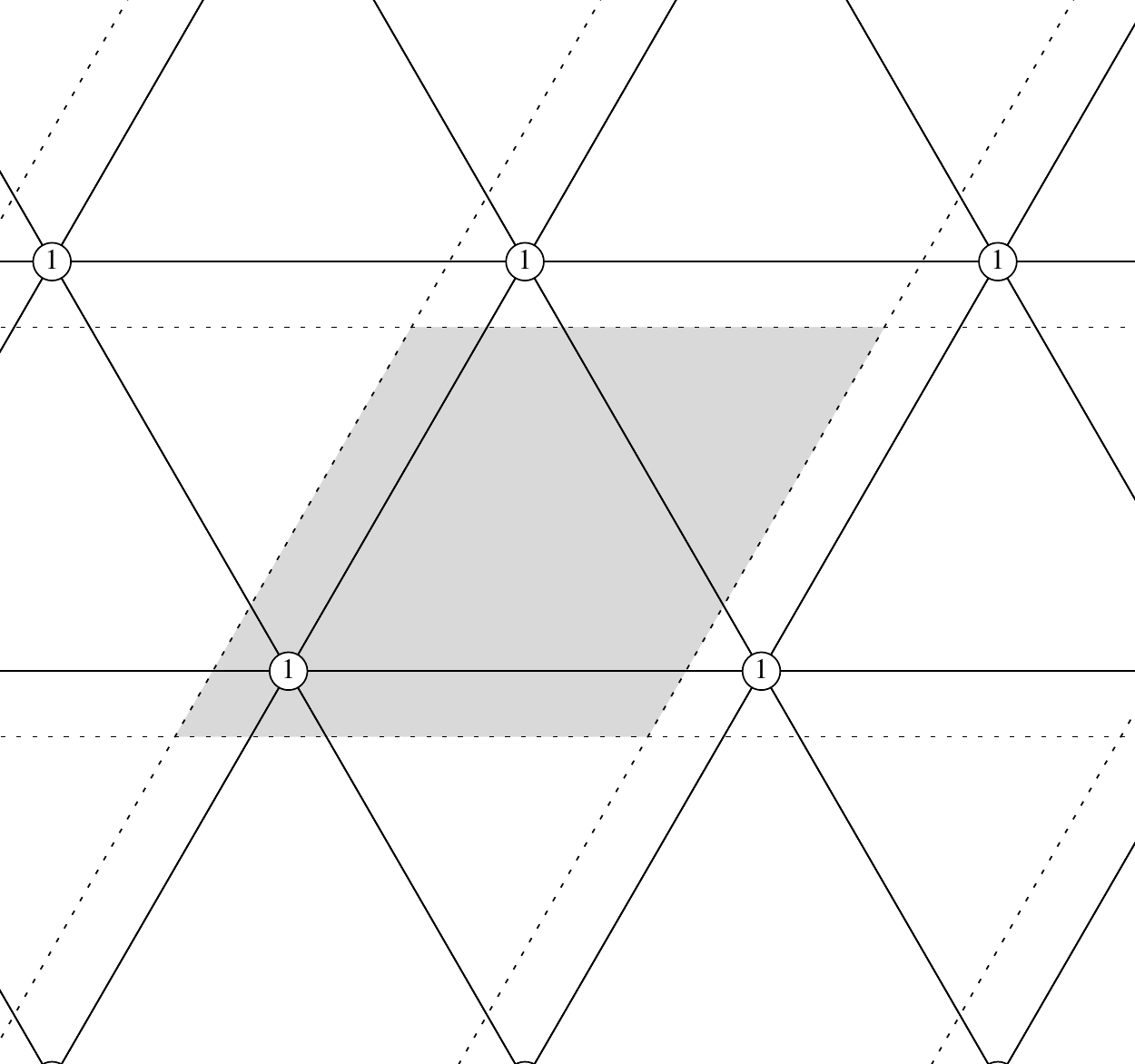}\hfill
\includegraphics[width=0.35\textwidth]{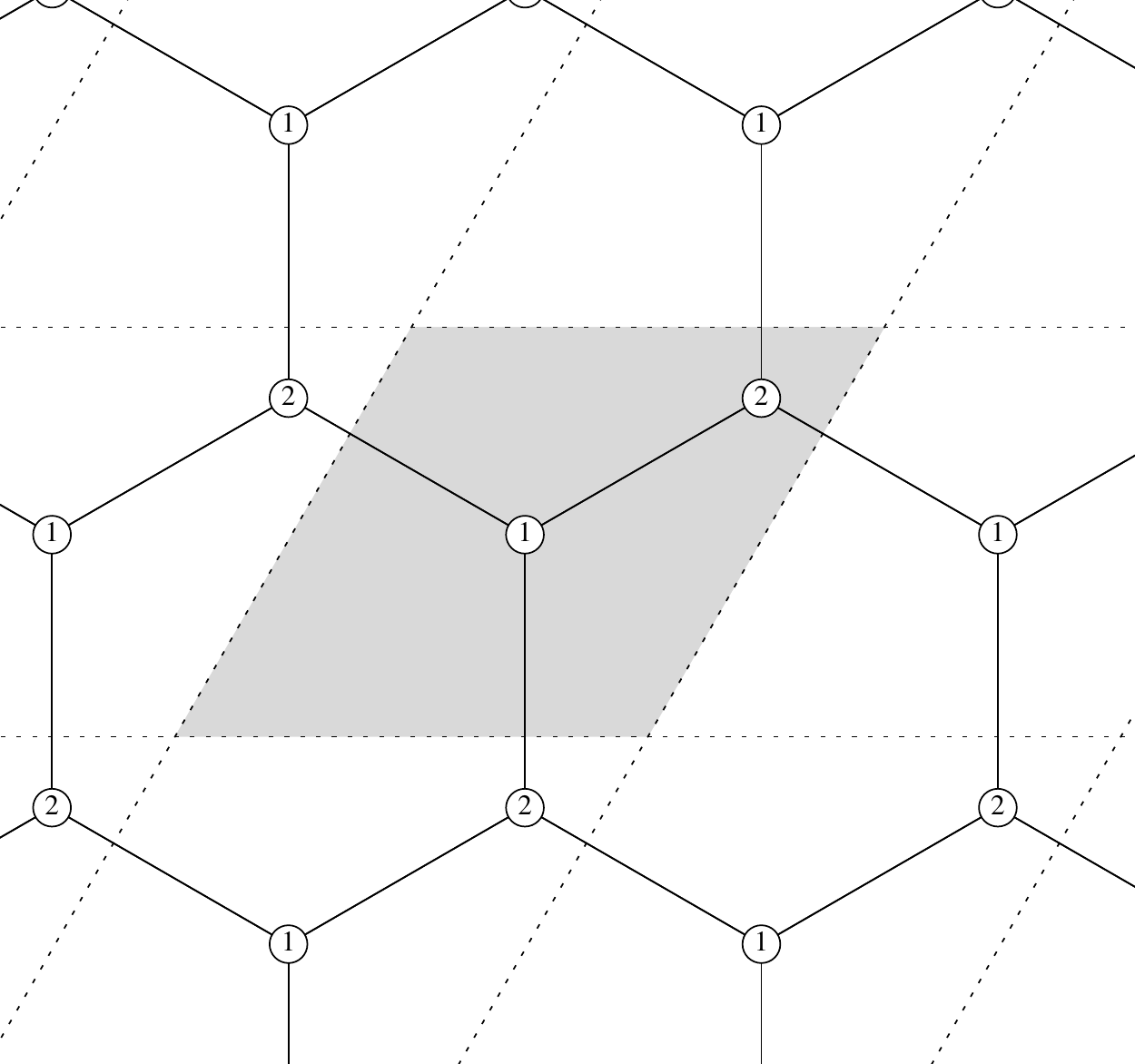}
\hfill}
\caption{Triangular and honeycomb lattices.}\label{triangular}
\end{figure}

\begin{figure}[t]
\centering{\hfill
\includegraphics[width=0.35\textwidth]{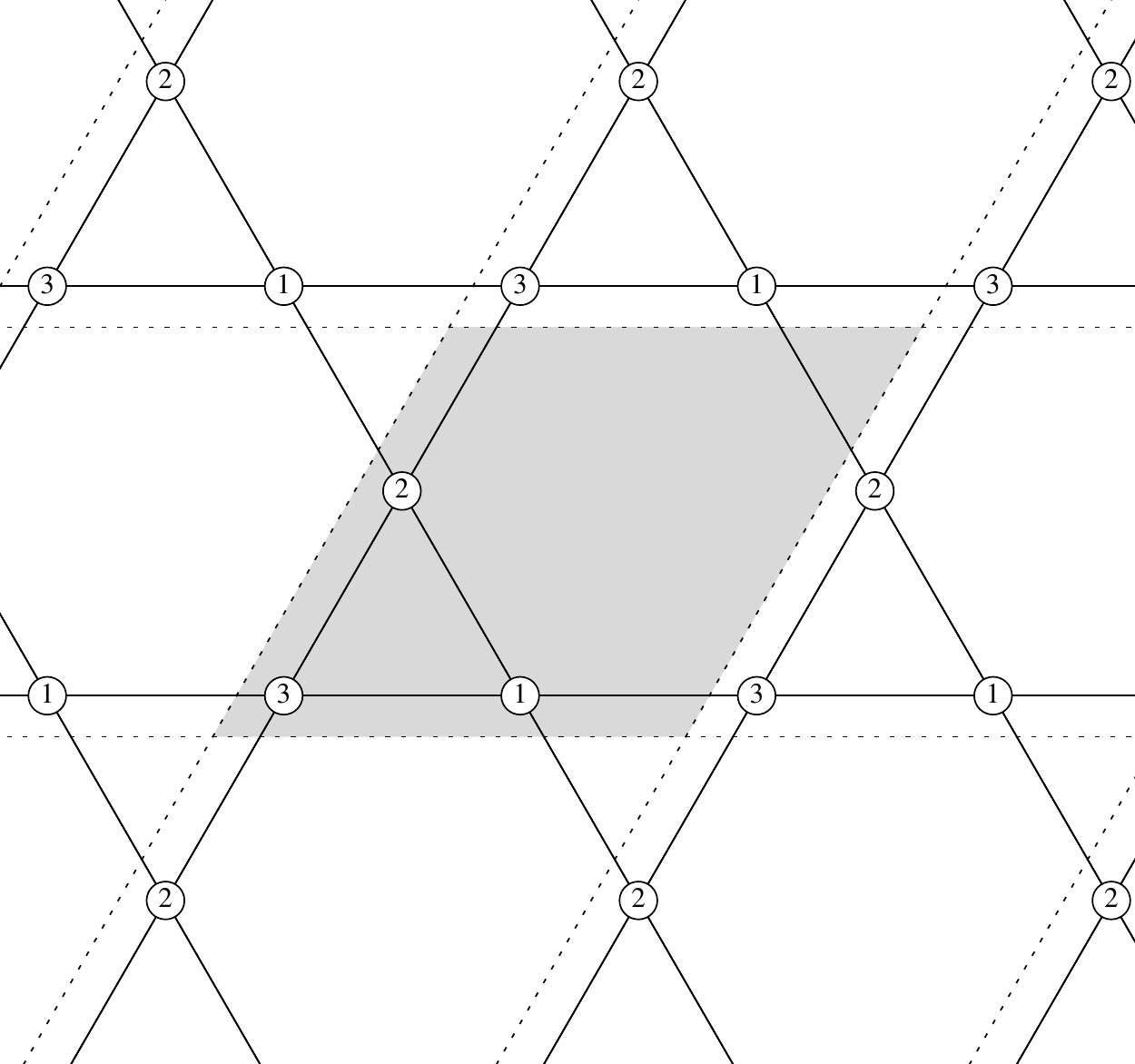}\hfill
\includegraphics[width=0.35\textwidth]{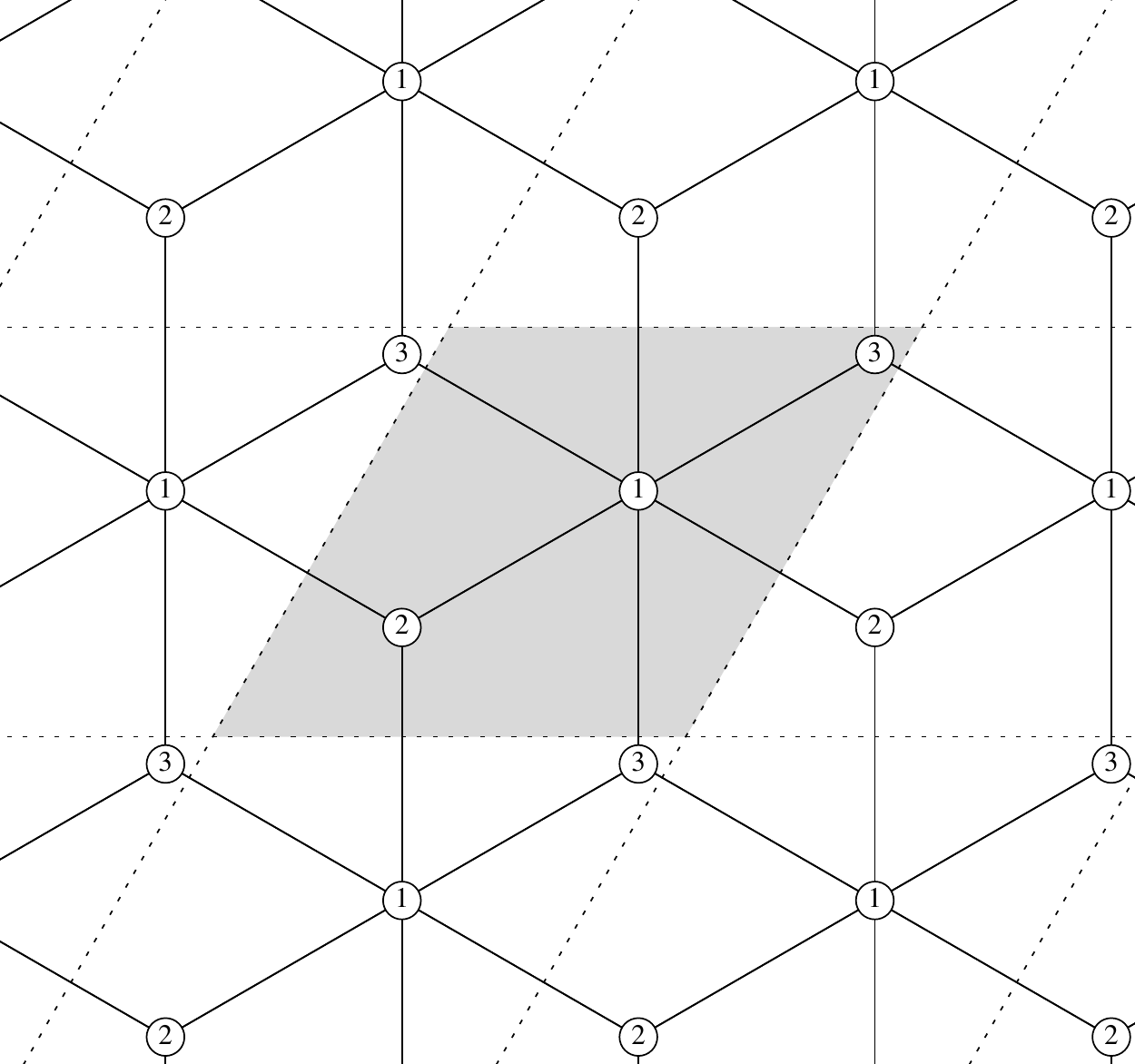}\hfill}
\caption{Kagom\'e and dice lattices.}\label{kagome}
\end{figure}

For the dice lattice (Figure~\ref{kagome}), the dual is the kagom\'e lattice, for which the potential kernel across edges is $1/4$ for each edge, and we have
\[ \tau\text{\tiny(dice)} = \frac{1}{4} \times \frac{1}{4} + \left(\frac14-\frac14\right)^2 = \frac{1}{16}\,.\]

The square, triangular, honeycomb, and dice lattices are all
sufficiently symmetric that $\tau=1/(\delta^*)^2$.
For general unweighted finite planar graphs,
\eqref{tau} together with the bound $a^2-ab+b^2 = \frac{1}{4}(a+b)^2 + \frac{3}{4}(a-b)^2 \geq \frac{1}{4}(a+b)^2$ imply
\[
  \tau\geq\sum_{e^*} \frac{\frac{1}{4}\Pr[e^*\in T^*]^2}{|E^*|} \geq \frac{|E^*|((|V^*|-1)/|E^*|)^2}{4|E^*|}\,,
\]
and hence $\tau\geq 1/(\delta^*)^2$ and also $\lambda\leq 2\delta^*$
for unweighted $\Z^2$-periodic planar lattices.
These bounds are all equalities for the square, honeycomb, triangular, and dice lattices,
and appear to be fairly good when each face of the lattice has the same number of sides.

\vspace{9pt}

For the kagom\'e lattice, the dual is the dice lattice (see Figure~\ref{kagome}), for which there are two types of vertices (degree-$3$ and degree-$6$), but only one type of edge.  By symmetry, for each edge~$(u,v)$, the potential kernel $A_{u,v}$ is $1/\operatorname{degree}(u)$.  Hence
\[ \tau\text{\tiny(kagom\'e)} = \frac{1}{6} \times \frac{1}{3} + \left(\frac16-\frac13\right)^2 = \frac{1}{12}\,.\]

\begin{figure}[t]
\centering{\hfill
\includegraphics[width=0.35\textwidth]{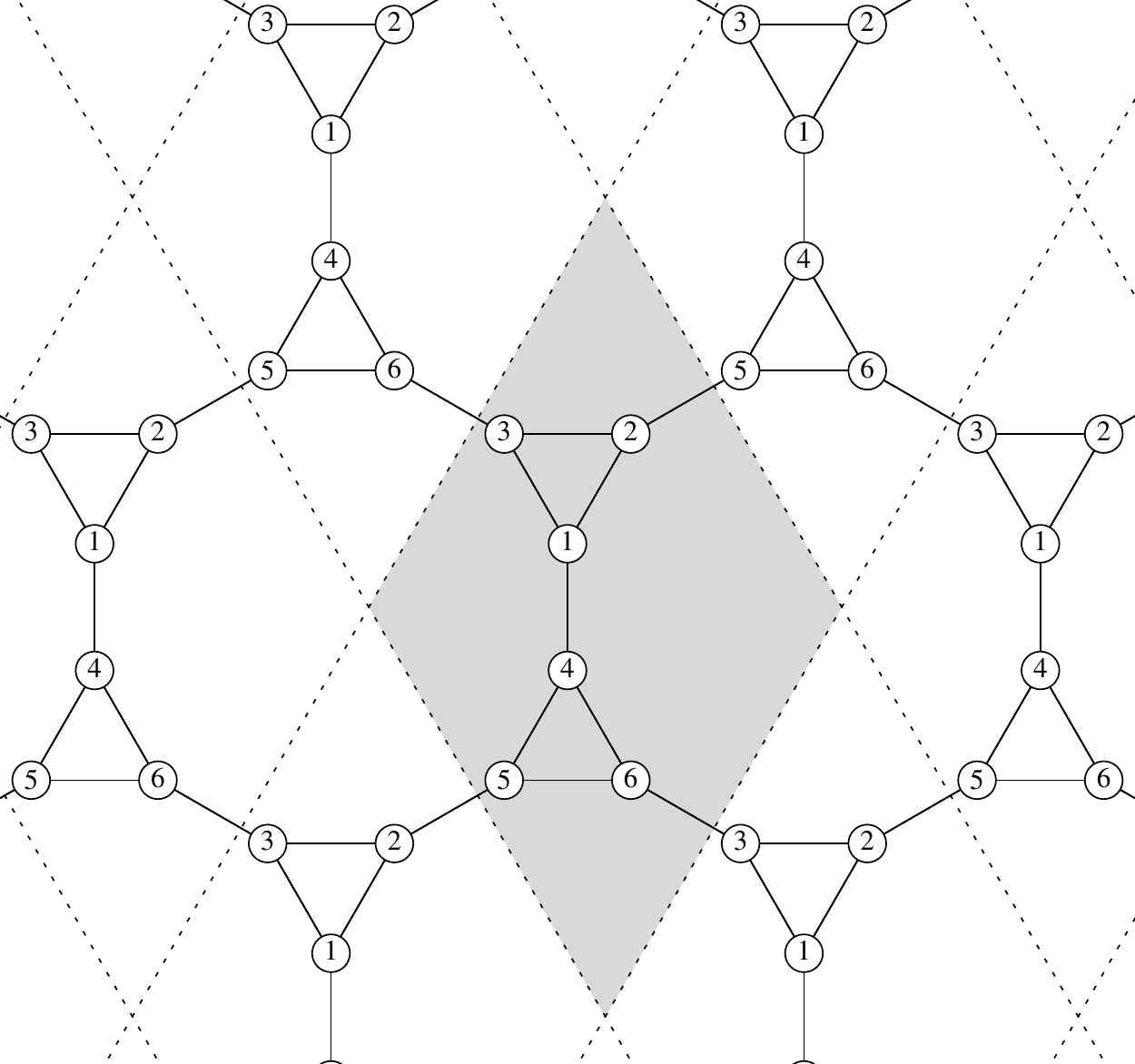}\hfill
\includegraphics[width=0.35\textwidth]{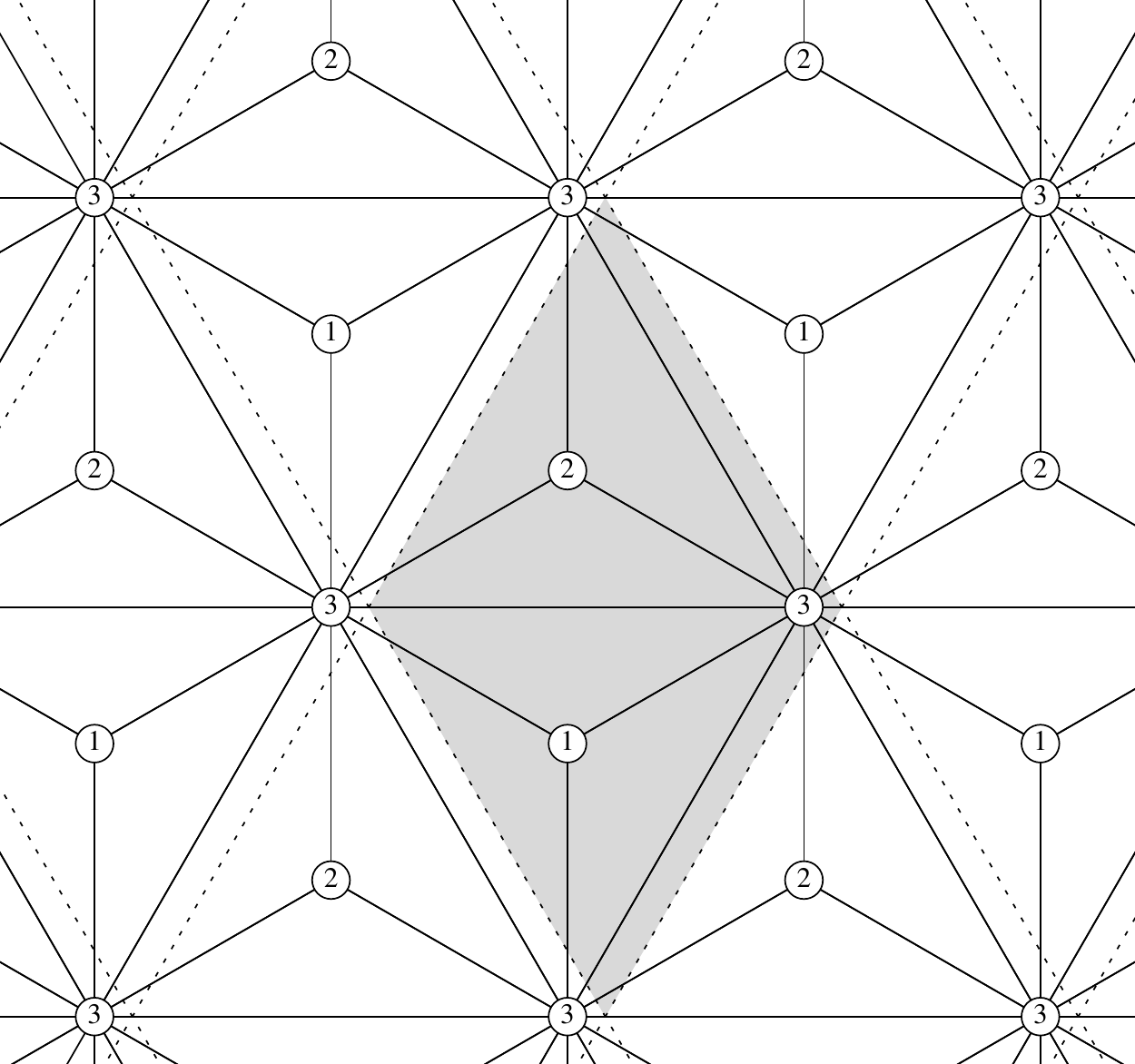}\hfill}
\caption{Fisher lattice and its dual, the triakis triangular lattice.}\label{fisher}
\end{figure}

The next pair of lattices that we consider are the Fisher lattice (i.e., the truncated
hexagonal lattice) and its dual, the triakis triangular lattice (see
Figure~\ref{fisher}).  There are several ways to determine the
potential kernel for adjacent vertices of these lattices; we describe
a way which essentially only uses symmetry.  In the triakis triangular
lattice, for the degree-3 vertices the potential kernel is of course
$1/3$, by symmetry.  Each degree-12 vertex $o$ is surrounded by 6 other
degree-12 vertices which are symmetric to one another; let $a$ denote
one such vertex.
It is also surrounded by 6 degree-3 vertices which are symmetric to one another;
let $b$ denote one such vertex.
The potential kernel $A_{u,v}$ is harmonic as a function of $v$ except at $u$.
Vertex~$b$ is surrounded by $o$ and two degree-12 neighbors of $o$.
By harmonicity and symmetry, $A_{o,b}=\frac{1}{3} A_{o,o} + \frac{1}{3} A_{o,a} + \frac{1}{3} A_{o,a}
= \frac23 A_{o,a}$.
Thus $1 = 6\times A_{o,a} + 6\times A_{o,b} = 10 A_{o,a}$,
so $A_{o,a}=\frac{1}{10}$ and $A_{o,b}=\frac{1}{15}$.
Next we use the following
relation~\eqref{primal-dual-kernel} between the potential kernels of a graph and
its dual:
\begin{equation}\label{primal-dual-kernel}
  w_{u,v}\big(A^{(s)}_{u,v}+A^{(s)}_{v,u}\big) + w_{u^*,v^*}\big(A^{(s^*)}_{u^*,v^*} + A^{(s^*)}_{v^*,u^*}\big) = 1 \quad\text{(planar undirected graphs)}.
\end{equation}
This follows from \eqref{uv-in-T} because either edge $(u,v)$ is in the tree, or its dual edge $(u^*,v^*)$
is in the dual tree.
Together with the symmetry in the Fisher lattice, we find
that the potential kernel along intertriangle edges is $\frac12
(1-\frac{1}{10}-\frac{1}{10}) = \frac{2}{5}$.  For the intratriangle
edges, we again use symmetry to find that the potential kernel is
$\frac{1}{2}(1-\frac{2}{5}) = \frac{3}{10}$.

For the triakis triangular lattice,
the 3--12 edges are twice as numerous as the 12--12 edges.  Thus,
for its dual the Fisher lattice, we obtain
\[ \tau\text{\tiny(Fisher)} = \frac{1}{3} \left[\frac{1}{10}\times\frac{1}{10} + \left(\frac{1}{10}-\frac{1}{10}\right)^2\right] + \frac{2}{3}\left[\frac{1}{3}\times\frac{1}{15} + \left(\frac{1}{3}-\frac{1}{15}\right)^2\right] = \frac{59}{900}\,.\]

For the Fisher lattice, the intratriangle edges are twice as numerous as the
intertriangle edges.  Thus, for its dual the triakis triangular lattice, we obtain
\[ \tau\big(\parbox{\widthof{\tiny triangular}}{\centering\tiny triakis\\[0pt]triangular}\big)
 = \frac{1}{3} \left[ \frac{4}{10}\times\frac{4}{10}\right] + \frac{2}{3} \left[ \frac{3}{10}\times\frac{3}{10} \right] = \frac{17}{150}\,.\]

There are two types of edges in these lattices, and it is natural to give them different edge weights.
If we give weight $\beta$ to the edges connecting two degree-12 vertices of the triakis triangular lattice,
and weight~$1$ to the other edges, then the same symmetry argument can be used to compute the
potential kernel.  This leads to
\begin{align*}
 \tau\big(\parbox{\widthof{\tiny triangular}}{\centering\tiny triakis\\[-0pt]triangular}\big)
 &= \frac{2\beta^3+8\beta^2+6\beta+1}{2(\beta+2)(3\beta+2)^2}\,,\text{ and} \\[6pt]
 \rho\big(\parbox{\widthof{\tiny triangular}}{\centering\tiny triakis\\[-0pt]triangular}\big)
 &= \frac{(\beta+1)(5\beta^2+11\beta+5)}{2(\beta+2)(3\beta+2)^2} \,.
\end{align*}
With $\beta=1$ we recover the parameters for the unweighted triakis
triangular lattice.  Since a random walk on the weighted triakis triangular lattice converges to
a random walk on the dice lattice as $\beta\to0$, and to a random walk on the triangular
lattice as $\beta\to\infty$, we can recover $\tau$ and $\rho$ for
the dice and triangular lattices by taking these limits.
There are similar rational function formulas for $\tau$ and $\rho$
on the weighted Fisher lattice, but it is not as simple to recover the parameters for the
honeycomb and kagom\'e lattices.

\begin{figure}[t!]
\centering{\hfill
\includegraphics[width=0.35\textwidth]{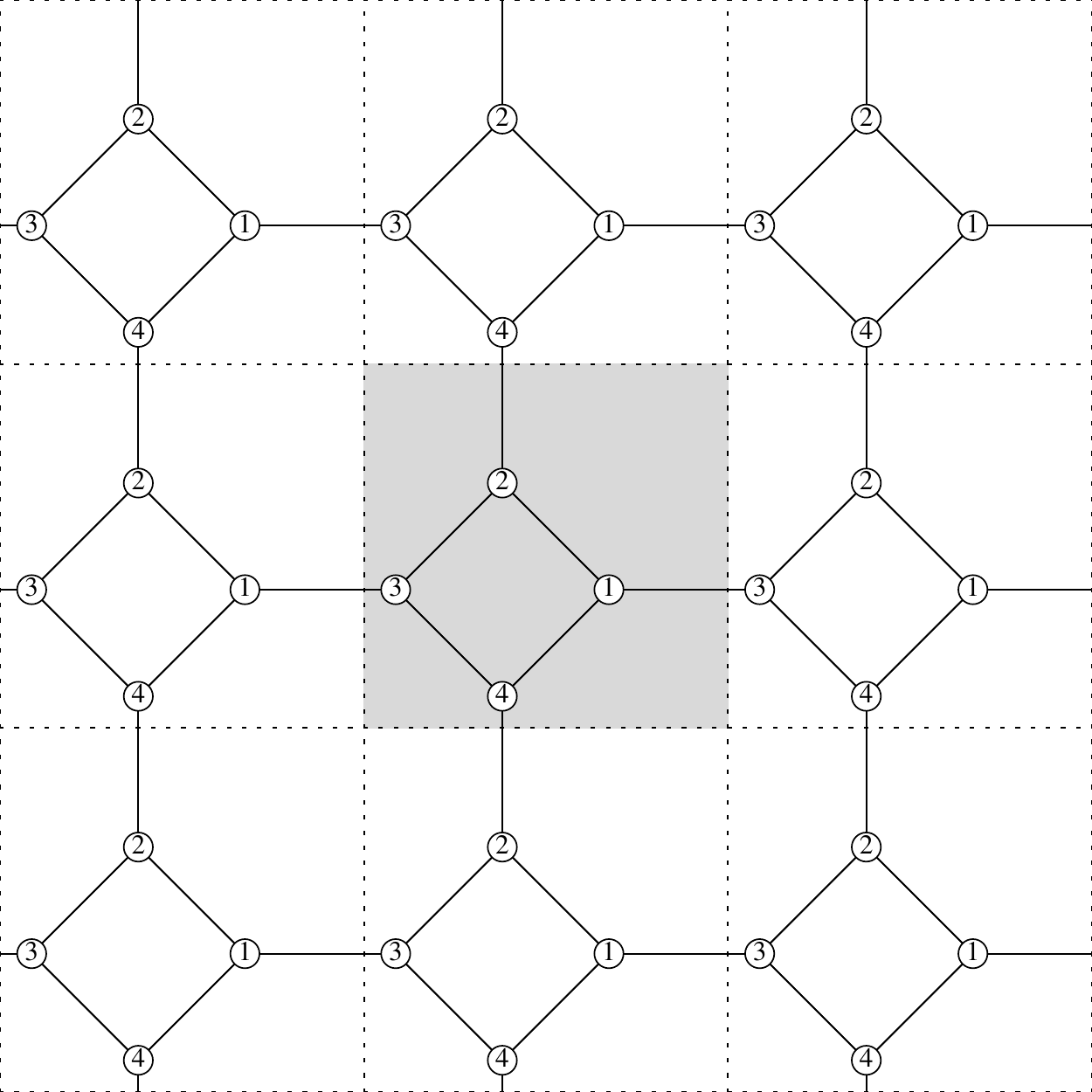}\hfill
\includegraphics[width=0.35\textwidth]{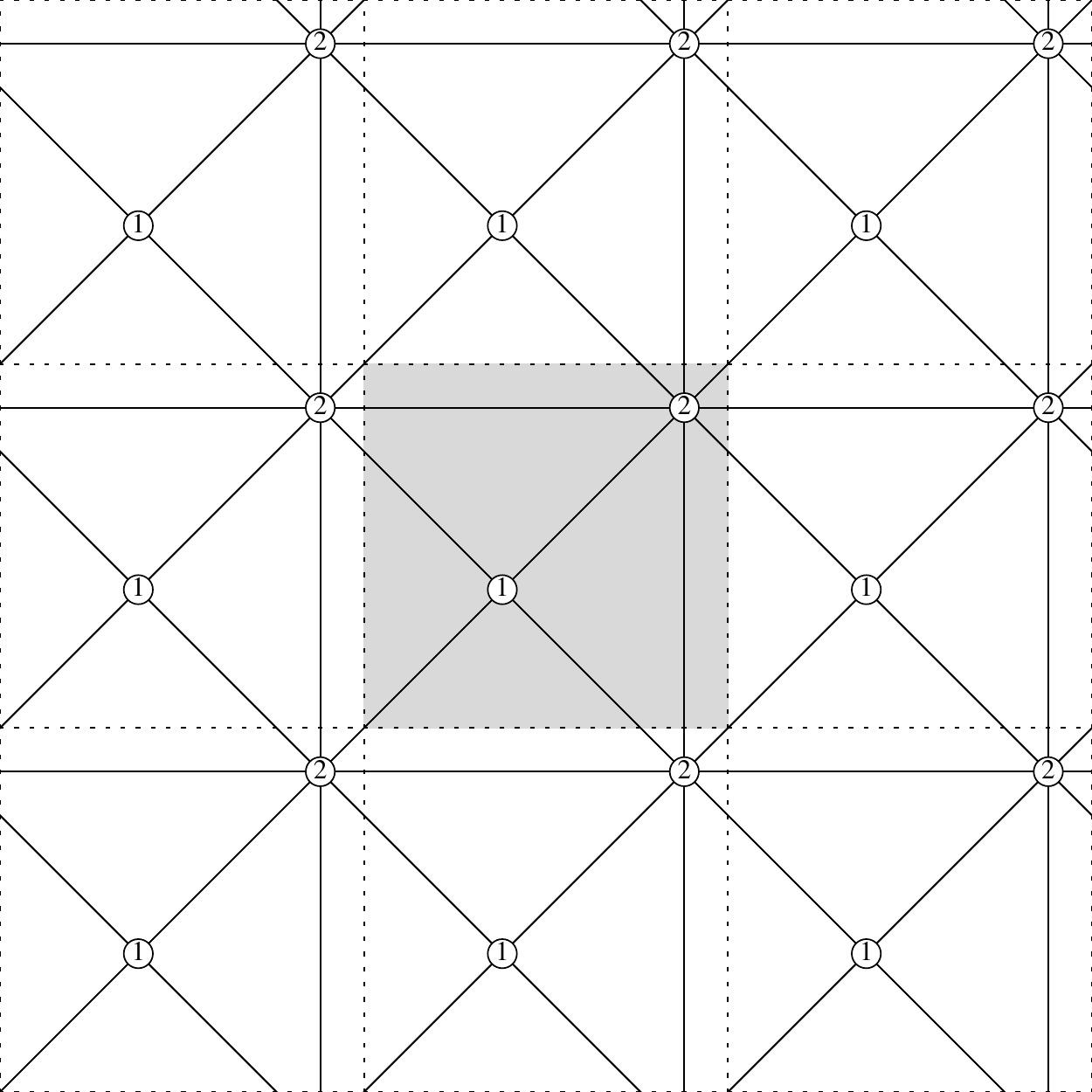}\hfill}
\caption{Square-octagon lattice and its dual, the tetrakis square lattice.}\label{square-octagon}
\end{figure}

The last pair of periodic lattices we consider are the square-octagon lattice (i.e., the truncated square lattice), and its dual the tetrakis square lattice (see Figure~\ref{square-octagon}).
Since there are two edge types, we give the edges connecting degree-8 vertices
in the tetrakis square lattice an edge weight of $\beta$, or equivalently, we give the edges
between octagons in the square-octagon lattice weight $1/\beta$.
In the tetrakis square lattice, the potential kernel between a
degree-4 vertex and one of its neighbors is $\frac{1}{4}$ by symmetry.
Let $\alpha=\alpha(\beta)$ denote the potential kernel between two adjacent degree-8
vertices.  Then the potential kernel from a degree-8 to a degree-4
vertex is $\frac{1}{4}-\alpha\beta$.  For the square-octagon lattice, we use
Equation~\eqref{primal-dual-kernel} relating the potential kernel of a graph to that of its
dual together with the bilateral symmetry of the edges to deduce that
for the intersquare edges the potential kernel is
$\frac12\beta(1-\alpha\beta-\alpha\beta)=\frac{1}{2}\beta-\alpha\beta^2$, from which it follows that the
potential kernel for the edges in the squares is $\frac{1}{4}+\alpha\beta/2$.
Substituting these values into Equations~\eqref{tau}, \eqref{eq:rho-tau}, and \eqref{eq:r-e-in-T},
and using the fact that the unweighted edges are twice as numerous as the weighted ones, we obtain
\begin{align*}
 \tau\big(\parbox{\widthof{\tiny tetrakis}}{\centering\tiny tetrakis\\[1pt]square}\big)
 &= \frac{\tfrac14\beta+\tfrac18}{\beta+2}
  - \frac{\beta-\tfrac12}{\beta+2}\beta\alpha
  + \frac{\beta+\tfrac12}{\beta+2}\beta^2\alpha^2\,,\\[6pt]
 \rho\big(\parbox{\widthof{\tiny tetrakis}}{\centering\tiny tetrakis\\[1pt]square}\big)
 &= \frac{\tfrac14\beta+\tfrac58}{\beta+2}
  - \frac{\tfrac32}{\beta+2}\beta\alpha
  + \frac{\beta+\tfrac12}{\beta+2}\beta^2\alpha^2\,,
\end{align*}
and
\begin{align*}
\tau\big(\parbox{\widthof{\tiny octagon}}{\centering\tiny square-\\[0pt]octagon}\big)
 &= \frac{\tfrac18\beta}{1+2\beta} - \frac{\tfrac12\beta}{1+2\beta} \beta\alpha + \beta^2\alpha^2\,,\\[6pt]
\rho\big(\parbox{\widthof{\tiny octagon}}{\centering\tiny square-\\[0pt]octagon}\big)
 &= \frac{\tfrac58\beta+\tfrac12}{1+2\beta} + \frac{\tfrac12\beta-1}{1+2\beta} \beta\alpha + \beta^2\alpha^2\,.
\end{align*}

Symmetry is not enough to determine the value of $\alpha$, but we can
use a method that is applicable to any periodic graph \cite{Spitzer}.
We write the Laplacian in Fourier space as a matrix
$\widehat{\Delta}(z,w)$ indexed by the vertices of a fundamental
domain.  An edge connecting a vertex of type~$i$ to a vertex of
type~$j$ in the fundamental domain $x$ units in the $z$-direction and
$y$ units in the $w$-direction contributes $1$ to $\Delta_{i,i}$ and
$-z^x w^y$ to~$\Delta_{i,j}$.  The tetrakis square lattice's
fundamental domain has two vertices, and in this case we have
\[
\widehat{\Delta}(z,w) =
\begin{bmatrix}
 4 & -\frac{1}{w z}{-}\frac{1}{w}{-}\frac{1}{z}{-}1 \\[4pt]
 - w z {-}w {-}z {-}1 & 4{+}4\beta{-}\beta w{-}\frac{\beta}{w}{-}\beta z{-}\frac{\beta}{z}
\end{bmatrix}\,.
\]
The Green's function is given by the inverse Fourier transform of
$\widehat{G}(z,w)=\widehat{\Delta}^{-1}(z,w)$, that is a double
integral over the circle of a rational function in $z$ and $w$.  The
potential kernel $A_{u,v}$ between a vertex $u$ of type $i$ and a
vertex $v$ of type $j$ in a fundamental domain shifted by $x$ units in
the $z$-direction and $y$ units in the $w$-direction is
\[
  A_{u,v} = \oint\!\!\oint \left[\big(\widehat{\Delta}(z,w)^{-1}\big)_{i,i}-\frac{\big(\widehat{\Delta}(z,w)^{-1}\big)_{i,j}}{z^x w^y}\right] \frac{dz}{2\pi iz} \,\frac{dw}{2\pi iw}\,.
\]

The evaluation of $\alpha$ corresponds to the case $i=j=2$, and by
symmetry $(x,y)$ can be any of $(1,0)$, $(-1,0)$, $(0,1)$, or $(0,-1)$.
Also $\big(\widehat{\Delta}(z,w)^{-1}\big)_{2,2}=4/P(z,w)$, where
\[
P(z,w) = 16\beta+12 - (4\beta+2)(z+z^{-1}+w+w^{-1}) - (z+z^{-1})(w+w^{-1})\,.
\]
Using symmetry, we can average the integrals for $(x,y)=(0,1)$ and $(x,y)=(0,-1)$:
\[
 \alpha = 2 \oint\!\!\oint \frac{2-w-w^{-1}}{P(z,w)} \,\frac{dz}{2\pi iz} \,\frac{dw}{2\pi iw}\,.
\]
Now for $b<1$
\[
\oint \frac{1}{1+b(z+z^{-1})/2} \frac{dz}{2\pi i z} = \frac{1}{\sqrt{1-b^2}}\,,
\]
so
\[
 \alpha = 2 \oint \frac{2-w-w^{-1}}{\sqrt{(16\beta+12 - (4\beta+2)(w+w^{-1}))^2 - 4(4\beta+2+w+w^{-1})^2}} \,\frac{dw}{2\pi iw}\,,
\]
which after changing variables letting $w=e^{2\pi i t}$ and simplifying gives
\[
\alpha = \frac{1}{2\sqrt{\beta^2+\beta}} \int_0^1 \sqrt{\frac{1 - \cos(2 \pi t)}{(3+2/\beta) - \cos(2 \pi t)}} \,dt\,.
\]
Now $1-\cos(2\pi t) = 2\sin^2(\pi t)$.  Let $u=\sin(\pi t)$.  Then $dt=du/(\pi\sqrt{1-u^2})$, so
\[
\alpha = \frac{1}{\pi\sqrt{\beta^2+\beta}} \int_0^1 \sqrt{\frac{u^2}{(1+1/\beta) +u^2}} \,\frac{du}{\sqrt{1-u^2}}\,.
\]
Let $(2+1/\beta)v^2=1-u^2$,
so $(2+1/\beta)v\,dv=-u\,du=\sqrt{2+1/\beta}\sqrt{1-u^2}\,dv$, and
\[
\alpha = \frac{1}{\pi\sqrt{\beta^2+\beta}} \int\limits_{\frac{1}{\sqrt{2+1/\beta}}}^0 \frac{-dv}{\sqrt{1-v^2}} = \frac{\arcsin \big(1/\sqrt{2+1/\beta}\big)}{\pi\sqrt{\beta^2+\beta}}\,.
\]
The sine-doubling formula gives
\begin{align*}
 2 \arcsin \left(\frac{1}{\sqrt{2+1/\beta}}\right) &= \arcsin\left(2\sqrt{\frac{1}{2+1/\beta}}\sqrt{\frac{1+1/\beta}{2+1/\beta}}\right) \\ &= \arcsin \sqrt{1-\frac{1/\beta^2}{(2+1/\beta)^2}} = \arcsec(2\beta+1)\,,
\end{align*}
so
\[
\alpha=\alpha(\beta) = \frac{\arcsec(2\beta+1)}{2\pi\sqrt{\beta^2+\beta}}\,.
\]

For every rational value of the edge weight $\beta$ except $\frac12$, the number $\alpha(\beta)$ is transcendental.
(The exceptional point is $\alpha(\frac12)=1/\sqrt{27}$.)  This follows from the Gelfond-Schneider theorem (see \cite[Chapt.~10]{Niven}) together with the result that for rational $r$, the only rational values of $\cos(r\pi)$ are $0,\pm1,\pm\frac12$ (see \cite[Cor.~3.12]{Niven}).  Consequently the sandpile density (and the parameters $\tau$ and $\rho$) are transcendental for $\beta\in\Q^+\setminus\{\frac12\}$.

\section{Open problems}\label{openproblems}

On the infinite branching tree with degree $\delta$ (and wired
boundary conditions), the sandpile density is $\delta/2$
\cite{dhar-majumdar:bethe}.  It would be interesting if the (wired)
sandpile density could be computed for other planar hyperbolic
lattices.

\medskip

The variance in the amount of sand of a recurrent sandpile
configuration is the variance in its level.  Expressing this in terms
of the binomial moments and using~\eqref{T'(1,y)} and planar duality, the variance is
\begin{equation}\label{varianceformula}
2 \frac{F_3(\G^*)}{F_1(\G^*)} + \frac{F_2(\G^*)}{F_1(\G^*)} - \left(\frac{F_2(\G^*)}{F_1(\G^*)}\right)^2,
\end{equation}
where the forest formula \cite{liu-chow,myrvold} gives
\begin{equation} \label{eq:F_k}
  \frac{F_k(\G)}{F_1(\G)}
 =
  \sum_{h=0}^{k-1} (-1)^h \!\!\!\!\!\!
  \sum_{\substack{\{u_1\sim v_1,\dots,u_h\sim v_h\}\subset E\\ \{x_1,\dots,x_{k-1-h}\}\subset V}}
  \!\!\!
  w_{u_1,v_1}\cdots w_{u_h,v_h}
 \det G_{u_1,v_1,\dots,u_{h},v_h,x_1,\dots,x_{k-1-h}}^{u_1,v_1,\dots,u_{h},v_h,x_1,\dots,x_{k-1-h}}
\end{equation}
where the inner sum is over sets of $h$ undirected edges and sets of $k-1-h$ vertices.
For the square grid, the variance in the amount of sand in a random
recurrent sandpile configuration of an $n\times n$ box in $\Z^2$
appears to be $(0.14386408549334\dots+o(1)) \times n^2$.
Is there a closed-form expression for this asymptotic variance?
Is the total amount of sand distributed according to a Gaussian?

\medskip

For non-planar graphs, what is
the complexity of counting spanning unicyclic subgraphs?
Is it \#P-hard?
Is it polynomial time solvable?  Is there a good formula which can be
used to find the sandpile density?

\section*{Acknowledgments}

We thank Igor Pak for bringing~\cite{myrvold} to our attention,
Henry Cohn for explaining to us why $\arcsec(p/q)/\pi$ is (usually) transcendental,
and the referees and editor for their suggestions.
A.K. acknowledges the hospitality of Microsoft Research Redmond where this work started.

\newcommand{\MRhref}[2]{\href{http://www.ams.org/mathscinet-getitem?mr=#1}{MR#1}}
\def\@rst #1 #2other{#1}
\newcommand\MR[1]{\relax\ifhmode\unskip\spacefactor3000 \space\fi
  \MRhref{\expandafter\@rst #1 other}{#1}}

\newcommand{\arXiv}[1]{\url{http://arxiv.org/abs/#1}}

\renewcommand\MR[1]{}

\providecommand{\bysame}{\leavevmode\hbox to3em{\hrulefill}\thinspace}
\providecommand{\MR}{\relax\ifhmode\unskip\space\fi MR }
\providecommand{\MRhref}[2]{%
  \href{http://www.ams.org/mathscinet-getitem?mr=#1}{#2}
}
\providecommand{\href}[2]{#2}


\label{LastPage}
\thispagestyle{empty}
\enlargethispage{40pt}

\begin{biog}
\item[Adrien Kassel] studied mathematics at \'Ecole Normale Sup\'erieure in Paris from 2006 to 2010. He obtained a doctoral degree in mathematics from Universit\'e Paris Sud in 2013.
\begin{affil}
ETH Z\"urich, Departement Mathematik, R\"amistrasse 101, 8092 Z\"urich, Switzerland\\
\texttt{adrien.kassel@math.ethz.ch}
\end{affil}

\item[David B.~\!Wilson] graduated from MIT in 1991 with degrees in electrical engineering, math, and computer science, and graduated again from MIT with a Ph.D.\ in math in 1996.  He specializes in probability and algorithms, and has a longtime interest in trees and sand.
\begin{affil}
Microsoft Research, One Microsoft Way, Redmond, WA 98052, U.S.A.\\
\url{http://dbwilson.com}
\end{affil}
\end{biog}
\vfill\eject

\end{document}